\documentclass[12pt]{article}
\usepackage[final]{graphicx}
\usepackage[dvips]{color}
\usepackage{amssymb}
\usepackage{graphicx}
\usepackage{amsmath}
\usepackage{amsfonts}
\usepackage{pifont}

\allowdisplaybreaks[1] 

\newcommand\e{{\varepsilon}}

\newcommand\D{{\Delta}} 

\newcommand\bu{{\bullet}} 
\newcommand\sq{{\square}}
\newcommand\ba{{\backslash}}
\newcommand\ci{{\circ}} 
\newcommand\bl{{\blacksquare}}

\newtheorem{cor}{Corollary}

\newtheorem{prp}{Proposition}

\def\comment#1{ }

\def\overset#1#2{{\buildrel #1\over#2}}

\begin{document}
\title{Structure of chambers cut out by Veronese arrangements of hyperplanes in the real projective spaces}
\author{Francois Ap\'ery, Bernard Morin and Masaaki Yoshida
}
\maketitle 
\section*{Introduction}
Consider several lines in the real projective plane. The lines divide the plane in various chambers. We are interested in the arrangement of the chambers. This problem is quite naive and simple. But if the cardinality $m$ of the lines are large, there are no way to control in general. If we assume that no three lines meet at a point (line arrangements {\bf in general position}), the situation does not improve much. In this paper we always assume this.  Let us observe when $m$ is small. Since we are mathematicians, we start from $m=0$: The projective plane itself. 
When $m=1$, what remains is the euclidean plane, of course. 
When $m=2$, the plane is divided in two di-angular chambers. 

When $m=3$, the lines cut out four triangles; if we think one of the lines is the line at infinity, then the remaining two can be considered as two axes dividing the euclidean plane in four parts. Another explanation: if we think that the three lines bound a triangle, then around the central triangle there is a triangle along each edge/vertex.

When $m=4$, four triangles and three quadrilaterals are arranged as follows; if we think one of the lines is the line at infinity, then the remaining three lines bound a triangle, and around this triangle there are a quadrilateral along each edge, and a triangle kissing at each vertex. Another explanation: put the four lines as the symbol $\#$, then around the central quadrilateral are a triangle along each edge, and a quadrilateral kissing at each  pair of opposite vertices. The situation is already not so trivial.

When $m=5$, there is a unique pentagon, and if we see the arrangement centered at this pentagon, the situation can be described simply: there are five triangles adjacent to five edges, and five quadrilaterals kissing at five vertices. If we think one of the lines is the line at infinity, the remaining four lines can not be like $\#$ since the five lines should be in general position; the picture does not look simple, but one can find a pentagon. Thanks to the central pentagon, we can understand well the arrangement.

When $m=6$, there are four types of arrangements. The simplest one has a hexagon surrounded by six triangles along the six edges; there are six quadrilaterals kissing at the six vertices, and three quadrilaterals away from the hexagon. Another one with icosahedral symmetry: If we identify the antipodal points of the dodeca-icosahedron projected from the center onto a sphere we get six lines in the projective plane: there are $20/2=10$ triangles and $12/2=6$ pentagons. We do not describe here the two other types.

For general $m$ there are so many types, and we can not find a way to control them, unless the $m$ lines bound an $m$-gon, from which we can see the arrangement. We can think this $m$-gon as the {\it center} of this arrangement, and we would like to study higher dimensional versions of this kind of arrangements with center, if such a chamber exists. Before going further, since we are mathematicians, we start from the very beginning: the $(-1)$-dimensional projective space is an empty set; well it is a bit difficult to find something interesting. The $0$-dimensional projective space is a point; this point would be the center, OK. Now we proceed to the projective line. $m$ points on the line divide the line into $m$ intervals. (In this case `in general position' means `distinct'.) Very easy, but there is no {\it center} in the arrangement! 

How about 3-dimensional case? We have $m$ planes in general position (no four planes meet at a point).  
$m=0$: The projective space itself. 
When $m=1$, what remains is the euclidean space, of course. 
When $m=2$, the space is divided into two di-hedral chambers. 
When $m=3$, the space is divided into four tri-hedral chambers. 

When $m=4$, the planes cut out eight tetrahedra; 
if we think one of the planes is the plane at infinity, then the remaining three can be considered as the coordinate planes dividing the euclidean space in eight parts. Another explanation: if we think that the four planes bound a tetrahedron, then around the central tetrahedron there is a tetrahedron adjacent to each face (kisses at the opposite vertex), and a tetrahedron touching along a pair of opposite edges: $1+4+6/2=8$.

When $m=5$, if we think one of the planes is the plane at infinity, then the remaining four bounds a tetrahedron. Around this central tetrahedron, there is a (triangular) prism adjacent to the face, a prism touching along an edge, and a tetrahedron kissing at a vertex. ($1T+4P+6P+4T.$) Since the central tetrahedron is bounded only by four planes, this can not be considered as a center of the arrangement. Though a prism is bounded by five planes, there are several prisms, and none can be considered as the center.

When $m=6$, we can not describe the arrangement in a few lines; we will see that there is no center of the arrangement. The last author made a study of this arrangement in \cite{CYY}, but he himself admits the insufficiency of the description. To make the scenery of this arrangement visible, we spend more than thirty pages. The chamber bounded by six faces which is not a cube plays an important role: this is bounded by two pentagons, two triangles and two quadrilaterals. This fundamental chamber seems to have no name yet. So we name it as a {\it dumpling} ({\it gyoza} in  \cite{CY, CYY}).

\par\medskip
We study arrangements of $m$ hyperplanes in the real projective $n$-space in general position. As we wrote already, there would be no interesting things in general. So we restrict ourselves to consider the {\it Veronese arrangements} (this will be defined in the text); when $n=2$, it means that $m$ lines bound an $m$-gon. In general such an arrangement can be characterized by the existence of the action of the cyclic group $\mathbb{Z}_m=\mathbb{Z}/m\mathbb{Z}$ of order $m$. Note the fact: if $m\le n+3$ then every arrangement is Veronese. 

When $n$ is even, there is a unique chamber which is stable under the group; this chamber will be called the {\it central} one. When $n=0$, it is the unique point, and when $n=2$, it is the $m$-gon. We are interested in the next case $n=4$, which we study in detail. In particular, when $(n,m)=(4,7)$, the central chamber is bounded by seven dumplings. This study gives a light to general even dimensional central chamber. Remember that central chambers will be  higher dimensional versions of a point and a pentagon (and an $m$-gon).

When $n$ is odd, other than the case $(n,m)=(3,6)$, we study higher dimensional versions of a dumpling. Note that the 1-dimensional dumpling should be an interval.
\tableofcontents
\section{Preliminaries}
We consider {\bf arrangements} of $m$ hyperplanes in the $n$-dimensional real projective space
$$\mathbb{P}^n:=\mathbb{R}^{n+1}-\{0\}\ {\rm modulo}\ \mathbb{R}^\times.$$
$\mathbb{P}^{-1}$ is empty, $\mathbb{P}^0$ is a singleton, $\mathbb{P}^1$ is called the projective line, $\mathbb{P}^2$ the projective plane, and $\mathbb{P}^n$ the projective $n$-space.
We always suppose the hyperplanes are {\it in general position} (i.e. no $n+1$ hyperplanes meet at a point). 
We often work in the $n$-dimensional sphere
$$\mathbb{S}^n:=\mathbb{R}^{n+1}-\{0\}\ {\rm modulo}\ \mathbb{R}_{>0},$$
especially when we treat inequalities. This is just the double cover of $\mathbb{P}^n$, which is obtained by identifying the anti-podal points of $\mathbb{S}^n$.
So, a hyperplane in $\mathbb{S}^n$ is nothing but a punctured vector hyperplane modulo $\mathbb{R}_{>0}$.
\subsection{Projective spaces}
Let us give some geometric idea of the projective spaces.
\begin{itemize}
\item Two distinct hyperplanes meet along a projective space two dimensional lower. Two distinct lines in $\mathbb{P}^2$ meet at a point. Two distinct points in $\mathbb{P}^1$ do not meet. 
\item If you think a hyperplane as the hyperplane at infinity, what remains is the usual euclidean space. A tubular neighborhood of a line -- the complement of a disc -- in $\mathbb{P}^2$ is a M\"obius strip.
\item If you keep walking along a straight line then you eventually come back along the same line from behind to the point you started.
\item Two parallel lines meet at a point at infinity.
\item $\mathbb{P}^1$ as well as $\mathbb{S}^1$ is just a circle.
\item $\mathbb{P}^2$ can not be embedded into $\mathbb{R}^3$ unless you permit a self-intersection.
\item Odd dimensional ones are orientable, and even dimensional ones are non-orientable.
\end{itemize}
\subsection{Hyperplane arrangements}
We consider arrangements of $m$ hyperplanes in $\mathbb{P}^n$. We always assume that an arrangement is {\it in general position}, which means no $n+1$ hyperplanes meet at a point. 
Let $x_0:x_1:\cdots:x_n$ be a system of homogeneous coordinates on $\mathbb{P}^n$. A hyperplane $H$ is defined by a linear equation
$$a_0x_0+a_1x_1+\cdots+a_nx_n=0,\quad (a_0,\dots,a_n)\not=(0,\dots,0).$$
By corresponding a hyperplane $H$  its coefficients $a_0:\cdots:a_n$ of the defining equation, we have an isomorphism between the set of hyperplanes 
and the set of points in the dual projective space (i.e. $a$-space). A hyperplane arrangement is {\it in general position} if and only if 
the corresponding point arrangement is in general position, which means no $n+1$ points are on a hyperplane. 
By definition an arrangement of $m$ hyperplanes $H_1,\dots,H_m$, is defined up to the action of the symmetry group on the indices $1,\dots,m$.
\subsubsection{Grassmann isomorphism}
An arrangement of $m(\ge n+2)$ hyperplanes
$$H_j: a_{0j}x_0+a_{1j}x_1+\cdots+a_{nj}x_n=0,\quad j=1,\dots,m$$
defines an $(n+1)\times m$-matrix $A=(a_{ij})$. We can regard $A$ as a matrix representing a linear map from a linear space of dimension $m$ to that of dimension $n+1$. 
Note that the arrangement is in general position if and only if no $(n+1)$-minor vanish. Choose a basis of the kernel of this map and arrange them vertically, 
and we get an $m\times(m-n-1)$-matrix $B$ such that $AB=0$. The matrix ${}^tB$ defines an arrangement of $m$ hyperplanes in $\mathbb{P}^{m-n-2}$. 
The choice of the bases is not unique; the ambiguity gives projective transformations of $\mathbb{P}^{m-n-2}$. 
No $(m-n-1)$-minor of $B$ vanish (linear (in)dependence of the first $n+1$ columns of $A$ implies that of the last $m-n-1$ lines of $B$).

\noindent Summing up, we get an isomorphism between the arrangements of $m$ hyperplanes in $\mathbb{P}^n$ and those in  $\mathbb{P}^{m-n-2}$.

\subsubsection{Arrangements when $m=n+2$}
It is well-known that any three distinct points on the projective line $\mathbb{P}^1=\mathbb{R}\cup\{\infty\}$ can be transformed projectively into $\{0,1,\infty\}$. 
In general, we want to prove that any two systems of $n+2$ points in general position in $\mathbb{P}^n$ can be transformed projectively to each other. 
Since Grassmann isomorphism doesn't make sense when $m-n-2=0$, we prove this fact using linear algebra.
\begin{prp}\label{projectivelyeq}
Any $n+2$ points in general position in $\mathbb{P}^n$ can be transformed projectively into the $n+2$ points:
$$1:0:\cdots:0,\quad0:1:0:\cdots:0,\quad\dots,\quad0:\cdots:0:1,\quad 1:\cdots:1.$$
\end{prp}
{\bf proof:}
Let $A$ be the corresponding $(n+1)\times (n+2)$-matrix. Multiplying a suitable matrix from the left, we can assume that $A$ is of the form
$(I_{n+1}\ a)$, where $a=^t(a_0,\dots,a_n)$. Since the arrangement is in general position, $a_j\not=0\ (j=0,\dots,n)$. 
We then multiplying diag$(1/a_0, \dots, 1/a_n)$ from the left.

\noindent Projective transformations still operate on these points as permutations of $n+2$ points. 
\subsubsection{Arrangements when $m=n+3$}
\begin{prp}The set of arrangements of $n+3$ hyperplanes in general position in $\mathbb{P}^n$ is connected. \end{prp}

\noindent We can use the Grassmann isomorphism (note that $m-n-2=1$). The arrangements of $m$ points in $\mathbb{P}^1$ is topologically unique, but not projectively of course. 
We give direct proof without using Grassmann isomorphism.
We prove the dual statement: The set of arrangements of $n+3$ points in $\mathbb{P}^n$ in general position is connected. 

\noindent{\bf proof:}
Put $n+2$ points as in Proposition \ref{projectivelyeq}. Where can we put the ($n+3$)-th point so that the $n+3$ points are in general position, that is, 
no $n+1$ points are on a hyperplane? $n$ points out of these $n+2$ points span hyperplanes defined by:
$$x_i=0,\quad x_j=x_k\quad (i,j,k=1,\dots,n+1,\ j\not=k).$$
These hyperplanes are places which are not allowed to put the ($n+3$)-th point. These hyperplanes divide the space $\mathbb{P}^n$ into simplices (if non-empty) defined by
$$x_{i_1}<x_{i_2}<\cdots<x_{i_{n+1}},\quad \{i_1,\dots,i_{n+1}\}\subset
\{0,1,\dots,n+1\},$$
where $x_0=0,x_{n+1}=1$. The symmetric group on $n+2$ letters acts transitively on these simplices. This completes the proof.

\noindent We can rephrase these propositions as:
\begin{cor} If $m\le n+2$ there is only one arrangement up to linear transformations. If $m=n+3$ there is only one arrangement up to continuous move keeping the intersection pattern.
\end{cor}
\subsubsection{Curves of degree $n$}
Consider a curve $x(t)$ in $\mathbb{P}^n$ given by
$$x_0=x_0(t),\quad x_1=x_1(t),\quad\dots,\quad x_n=x_n(t),\qquad t\in\mathbb{R}.$$
At $t=\tau$, if there is a hyperplane
$$y_0(\tau)x_0+y_1(\tau)x_1+\cdots+y_n(\tau)x_n=0$$
passing through $x(\tau)$ such that the derived vectors $x'(\tau), x''(\tau),\dots,x^{(n-1)}$ lie on it, that is,
$$X(\tau)
\left(\begin{array}{c}y_0(\tau)\\y_1(\tau)\\\vdots \\ y_n(\tau)\end{array}\right)=0,\quad
X(\tau)=\left(\begin{array}{cccc}
x_0(\tau)&x_1(\tau)&\cdots&x_n(\tau)\\
x'_0(\tau)&x'_1(\tau)&\cdots&x'_n(\tau)\\
\cdots&&&\\
x^{(n-1)}_0(\tau)&x^{(n-1)}_1(\tau)&\cdots&x^{(n-1)}_n(\tau)
\end{array}\right),
$$
 this hyperplane is called an {\it osculating hyperplane} of the curve at $x(\tau)$. If the rank of $X(\tau)$ is $n$, there is a unique osculating hyperplane at $x(\tau)$. The correspondence
$$\tau\longmapsto y_0(\tau):y_1(\tau):\cdots:y_n(\tau)$$
defines a curve in the dual projective space (i.e. $y$-space). Notice that
(since $(\sum x_ky_k)'=\sum x'_ky_k+\sum x'_ky_k, (\sum x'_ky_k)'=\sum x''_ky_k+\sum x'_ky'_k,$ etc) we have
$$Y(\tau)
\left(\begin{array}{c}x_0(\tau)\\\vdots \\ x_n(\tau)\end{array}\right)=0,\quad
Y(\tau)=\left(\begin{array}{ccc}
y'_0(\tau)&\cdots&y_n(\tau)\\
y_0(\tau)&\cdots&y'_n(\tau)\\
\cdots&&\\
y^{(n-1)}_0(\tau)&\cdots&y^{(n-1)}_n(\tau)
\end{array}\right),
$$
that is, the hyperplane $$x_0(\tau)y_0+\cdots+x_n(\tau)y_n=0$$in $y$-space is an osculating hyperplane of the curve $y(t)$ at $y(\tau)$. 
This curve $y(t)$ is called the {\it dual curve}. Notice also that if $x(t)$ is of degree $n$, so is $y(\tau)$. Thanks to this dual correspondence, in place of proving 
\begin{prp}
\label{pr1}
For any $n+3$ hyperplanes in general position in $\mathbb{P}^n$, there is a unique rational curve of degree $n$ osculating these hyperplanes.
\end{prp}
we prove
\begin{prp}For any $n+3$ points in general position in $\mathbb{P}^n$, there is a unique rational curve of degree $n$ 
passing through these points.\end{prp}
This is a generalization of a well-known fact that there is a unique conic passing through given five points in general position.

Without loss of generality, we put $n+3$ points as:
$$\begin{array}{cccccc}
    &x_1:&x_2&:\cdots:&x_n:&x_{n+1}\\
&&&&&\\
p_0=&1:&1&:\cdots:&1:&1,\\
p_1=&1:&0&:\cdots:&0:&0,\\
p_2=&0:&1&:\cdots:&0:&0,\\
\vdots&&&&&\\
p_n=&0:&0&:\cdots:&1:&0,\\
p_{n+1}=&0:&0&:\cdots:&0:&1,\\
p_{n+2}=&a_1:&a_2&:\cdots:&a_n:&a_{n+1},
\end{array}$$
where $0<a_1<\cdots<a_n<a_{n+1}$. Sorry, in this proof we use coordinates $x_1:\cdots:x_{n+1}$ instead of $x_0:\cdots:x_n$. We will find a curve
$$C:t\longmapsto x_1(t):\cdots:x_{n+1}(t),$$
such that $x_j(t)$ is a polynomial in $t$ of degree $n$, and that
$$ C(q_0)=p_0,\quad C(q_1)=p_1,\ \dots,\ C(q_{n+1})=p_{n+1},\quad C(r)=p_{n+2}.$$
If we normalize as $$q_0=\infty, \quad q_1=0, \quad q_2=1,$$
then the above condition is equivalent to the system of equations
$$\begin{array}{lll}
(x_1(r)=)&c(r-q_2)(r-q_3)(r-q_4)\cdots(r-q_{n+1})&=a_1,\\
(x_2(r)=)&c(r-q_1)(r-q_3)(r-q_4)\cdots(r-q_{n+1})&=a_2,\\
(x_3(r)=)&c(r-q_1)(r-q_2)(r-q_4)\cdots(r-q_{n+1})&=a_3,\\
\vdots&&\\
(x_{n+1}(r)=)&c(r-q_1)(r-q_2)(r-q_4)\cdots(r-q_{n})&=a_{n+1},
\end{array}$$
with $n+1$ unknowns $q_3,\dots,q_{n+1}, r$ and $c$.
From the first and the second equations, $r$ is solved, from the second and the third equation, $q_3$ is solved,..., and we obtain a unique set of solutions:
$$r=\frac{a_2}{a_2-a_1},\quad q_j=\frac{(a_j-a_1)a_2}{(a_2-a_1)a_j}\quad (j=3,\dots n+1)$$
We do not care the value of $c$.
Since 
$$q_3-1=\frac{(a_3-a_2)a_1}{(a_2-a_1)a_3},\quad q_j-q_i=\frac{a_1a_2}{a_2-a_1}\cdot\frac{a_j-a_i}{a_ja_i},\quad r-q_j=\frac{a_1a_2}{a_2-a_1}\cdot\frac1{a_j}$$
we have
$$q_1=0<q_2=1<q_3<\cdots<q_{n+1}<r.$$

\subsection{Number of chambers}
Let  $\#(n,m)$ be the number of chambers cut out from the projective $n$-space $\mathbb{P}^n$ by $m$ hyperplanes in general position. We have
$$\#(1,m)=m,\quad \#(n,1)=1.$$
Consider $m-1$ hyperplanes in $\mathbb{P}^n$. The $m$-th hyperplane intersects $m-1$ hyperplanes, which cut outs $\#(n-1,m-1)$ chambers of dimension $n-1$. Each $(n-1)$-dimensional chamber cut an $n$-dimensional chamber into two. Thus we have
$$ \#(n,m)=\#(n,m-1)+\#(n-1,m-1),$$
and so
$$\#(n,m)=1+\sum_{\ell=1}^{m-1}\#(n-1,\ell).$$
This can be readily solved:
$$\begin{array}{lll}
\#(n,m)&=\displaystyle{
{m\choose{n}}+{m\choose{n-2}}+\cdots+{m\choose2}+{m\choose0}},&n:{\rm even},\\[6mm]
\#(n,m)&=\displaystyle{{m\choose{n}}+{m\choose{n-2}}+\cdots+{m\choose3}+{m\choose1}},&n:{\rm odd},\\
  \end{array}$$
where ${m\choose n}=0$ if $m<n$. In particular, $\#(1,m)=m$,
$$\#(2,m)={m\choose2}+1,\quad \#(3,m)={m\choose3}+m,\quad\#(4,m)={m\choose4}+{m\choose2}+1.$$
\subsection{Rotating a polytope around a face}
\begin{figure}
  \begin{center}
  \includegraphics[width=90mm]{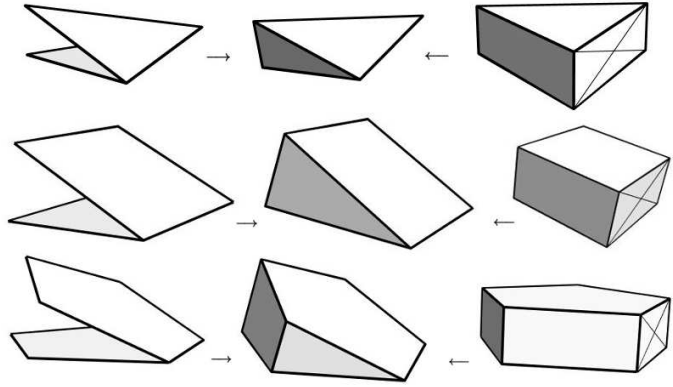}
  \end{center} 
  \caption{Dumplings $n=3, m=3,4,5$} 
  \label{dumpling}
\end{figure}

Let $P$ be an $n$-polytope in $\mathbb{R}^{n}$ canonically embedded in $\mathbb{R}^{n+1}$, and $\sigma$ one of its facets ($(n-1)$-dimensional faces).
The convex hull of $P$ and its image by a rotation in $\mathbb{R}^{n+1}$ of an angle less than $\pi$ and centered at $\sigma$ is an $(n+1)$-polytope, 
which will be denoted by $P^\sigma$, and will be called a polytope obtained from $P$ by {\it rotating} it around $\sigma$. It can be abstractly defined as
$$
P^\sigma=P\times[0,1]/\overset{\sigma}{\sim}, \quad (x,t)\ \overset{\sigma}{\sim}\ (x,0), \ x\in\sigma,\quad t\in[0,1],
$$
in other words it is obtained by crushing the facet $\sigma\times[0,1]$ of the product $P\times[0,1]$ by the projection (see Figure \ref{dumpling})
$$\sigma\times[0,1]\ni(x,t)\longmapsto x\in\sigma.$$ 
Note that the boundary of $P^\sigma$ consists of
$$P\times\{0\},\quad P\times\{1\},\quad (\tau\times[0,1])/\overset{\sigma}{\sim}\quad (\tau:\ {\rm face\ of\ }P,\ \tau\not=\sigma).$$

\noindent For example, rotating a pentagon around a side, we get a polytope bounded by two pentagons, two triangles, and two quadrilaterals. 
In other words, it is obtained by crushing a rectangular face of the pentagonal prism. 
As is shown in Figure \ref{dumpling} (right), it looks like a dumpling, so it will be called a {\bf dumpling}. 
This is called a  {\it gyoza} in \cite{CY, CYY}. In general, rotating an $m$-gon around a side, we get a polytope bounded by two $m$-gon's, 
two triangles, and $m-3$ quadrilaterals; this polytope will also be called a (3-dimensional $m$-) dumpling. 
A 3-dumpling is a tetrahedron, a 4-dumpling is a triangular prism.

\noindent High dimensional ones will be also called {\it dumplings}.

\section{Veronese arrangements}
A hyperplane arrangement $A=\{H_j\}_{j=1,\dots,m}$ in $\mathbb{S}^n$ (or $\mathbb{P}^n$) is said to be {\bf Veronese} if
under a suitable linear change of coordinates, $H_j$ is given by $f_j=0$, where
$$f_j=f(t_j,x)=x_0+t_jx_1+t_j^2x_2+\cdots+t_j^nx_n,\quad j=1,\dots,m.$$
Here
$t_1<t_2<\cdots<t_m$
are real numbers, and $(x_0,\dots,x_n)$ are coordinates on $\mathbb{R}^{n+1}$.
Note that by proposition \ref{pr1}, if $m\le n+3$, every arrangement of hyperplanes in general position is Veronese.

\begin{prp}
For the curve $V_n$ of degree $n$ defined by
$$\mathbb{R}\ni t\longmapsto x(t)=(-t)^n:n(-t)^{n-1}:\cdots:{n\choose k}(-t)^{n-k}:\cdots:1,$$
the osculating hyperplane at the point $x(t_j)$ is given by $f_j=0$.
\end{prp} The curve $V_n$ is often called a Veronese embedding of $\mathbb{P}^1$ into $\mathbb{S}^n$ (or $\mathbb{P}^n$). 
This is the reason for coining Veronese arrangements. 
\par\medskip
\noindent From here to the end of next section, we work in $\mathbb{S}^n$. The binomial theorem tells
$$f_j(x(t))=(t_j-t)^n;$$
so that if $n$ is even, 
$$
f_j(x(t))> 0\quad {\rm if}\quad t\not=t_j,
$$
if $n$ is odd,
$$
f_j(x(t))< 0\quad {\rm if}\quad t<t_j,\qquad f_j(x(t))>0\quad {\rm if}\quad t>t_j.
$$

\subsection{Chambers cut out by Veronese arrangements}
The closure of each connected component of the complement of a hyperplane arrangement is  called a {\bf chamber}. 
Let a Veronese arrangement be given as in the previous section.
Each chamber is given by a system of inequalities
$$\varepsilon_jf_j\ge0,\quad j=1,\dots,m, \quad(\varepsilon_j=\pm1)$$
which is often denoted by the sequence
$$\varepsilon_1,\ \dots,\ \varepsilon_m.$$
\comment{Remember we are working in $\mathbb{S}^n$. 
The facet (face of codimension 1) of the chamber on the hyperplane $H_k:f_k=0$ is given by
$$\varepsilon_jg_j\le0 \quad(j<k),\qquad \varepsilon_jg_j\ge0\quad (j>k);$$
this facet is often denoted by the sequence
$$-\varepsilon_1,\ \dots,\ -\varepsilon_{k-1},\ 0,\ \varepsilon_{k+1}.\ \dots,\ \varepsilon_m.$$}\par\noindent
{\bf Group action:} Since the Veronese arrangement is determined by a sequence $ t_1, t_2,\dots,t_m$ of points in $\mathbb{P}^1$ arranged in this order, the shift 
$$\iota: j\mapsto j+1 \quad{\rm mod}\ m$$ acts on the set of chambers by
$$\varepsilon_1,\ \dots,\ \varepsilon_m\longmapsto 
\left\{\begin{array}{ll}-\varepsilon_m,\ \varepsilon_1,\ \dots,\ \varepsilon_{m-1}\quad &m: {\rm \ odd,}\\
\varepsilon_m,\ \varepsilon_1,\ \dots,\ \varepsilon_{m-1}\quad &m: {\rm\ even.}
\end{array}\right. 
$$
The cyclic group generated by $\iota$ will be denoted by $\mathbb{Z}_m$. Note that if the points $t_1,t_2,\dots, t_m$ are so arranged that the transformation $t_j\to t_{j+1}$ is given by a projective one, then the corresponding action on $\mathbb{P}^n$ is also projective. Thus we have the well-defined $\mathbb{Z}_m$ action on the set of
chambers as above for any sequence $ t_1, t_2, \ldots , t_m \in \mathbb{P}^1 $.
\subsection{A specific feature of Veronese arrangements}
Let a Veronese arrangement $A=\{H_j\}$ of $m$ hyperplanes
$$H_j: f_j=f(t_j,x)=0,\quad \cdots<t_j<t_{j+1}<\cdots$$ 
in $\mathbb{P}^n$ be given as in the previous section. 
If we let $t_j$ tends near to $t_{j+1}$, the intersection pattern does not change, that is, 
there is no vertex (intersection point of $n$ hyperplanes in $A$) in the slit between the hyperplanes $H_j$ and $H_{j+1}$. 
More precisely, the hyperplanes $H_j$ and $H_{j+1}$ divide the space into two parts, and one of them does not contain any vertex. 

This fact gives the following information on the chambers cut out by $A$. If a chamber in the slit between the hyperplanes $H_j$ and $H_{j+1}$ 
does not touch the intersection $H_j\cap H_{j+1}$, then this chamber is the direct product of a chamber of the restricted arrangement $A_j$ and the unit interval. 
If a chamber in the slit touches the intersection of the two hyperplanes, 
then (since the hyperplanes are in general position) there is a chamber $P$ of the restricted arrangement $A_j$ 
with a facet  $\sigma$ included in $H_j\cap H_{j+1}$and the chamber in the slit is the polyhedron $P^\sigma$ obtained from $P$ by rotating around $\sigma$.
\par\medskip
A polyhedron (of dimension greater than 1) is said to be {\it irreducible} if it is not a direct product of two polyhedron nor is 
obtained from a lower one by  rotation described in \S 1.4. For example, triangles, rectangles, tetrahedra, prisms, cubes and dumplings are reducible, 
while pentagon is irreducible.
\par\medskip
A combinatorial study on Veronese arrangements is made in \cite{CY}. We do not use in this paper the result obtained there. To give a general idea, we just quote one of the results:
Consider a Veronese arrangement $A(m,n)$ of $m\ (\ge n+3)$ hyperplanes in $\mathbb{P}^n$ $(n\ge 2)$. 
\begin{itemize}
\item If $n$ is odd, then every chamber is a direct product of the unit interval and a chamber of $A(m-1,n-1)$, or is obtained from a chamber of $A(m-1,n-1)$ by rotating it around a facet.
\item If $n$ is even, then there is a unique irreducible chamber. Other chamber is obtained from $A(m-1,n-1)$ like the above.
\end{itemize}

\par\bigskip
In this paper, for a Veronese arrangement of even dimension, the unique irreducible chamber will be called the {\it central chamber}.

\section{Six planes in the 3-space}
\label{6planes}
 In this case there is no central chamber.
 The arrangement is Veronese. We are going to choose a particular arrangement so that
 the action of a group $G\cong\mathfrak{S}_3\times\mathbb{Z}_2$ of orientation 
 preserving projective transformations will be geometrically
 visible. 
 
 The elements of the group are given by matrices
 (up to scalar multiplication) with coefficients $0,1,-1$. The group $G$ is a subgroup of $PGL_+(4,\mathbb{Z})$. 
 Notice that $GL_+(4,\mathbb{Z})$ is not a group since inverse operation is not always defined.
 However it is defined up to scalar matrix multiplication, so that the quotient $PGL_+(4,\mathbb{Z})$ has a group structure.

 The group $G$ acts on the $26$ chambers with four orbits:
 $12$ prisms, $6$ tetrahedra, $6$ dumplings, $2$ cubes.
 The orbit of a prism together with a cube is a solid torus. The orbit of a dumpling is a 
 second solid torus linked with the first one. 
 The orbit of a tetrahedron is the complementary of the two solid tori.
 Actually, the union of one of these solid tori and the orbit of a tetrahedron is a solid torus as well.
 
\subsection{Simple observations}
\subsubsection{Cutting a tetrahedron and a prism}
Chambers in the plane bounded by five or less lines are
\begin{center}
diangle, triangle, quadrilateral,  and pentagon.
\end{center}
%
Chambers in the space bounded by four or less planes are dihedron, trihedron, tetrahedron; to see further we cut a tetrahedron by a plane. 
There are two ways to cut (see Figure \ref{tetracut}).
\begin{figure}[ht]
\begin{center} 
\include{tetracut}
\end{center}
\caption{Cutting a tetrahedron (tetracut)}
\label{tetracut}
\end{figure}
If we denote a tetrahedron by $T$, and a (triangular) prism by $P$, then the two cuttings can be presented as 
$$T\to T+P\quad({\rm triangle}),\quad P+P\quad({\rm quadrilateral}),$$
where $T+P$ (triangle) means $T$ and $P$ share a triangle.

We cut next a prism. There are five ways to cut (see Figure \ref{prismcut}). 
\begin{figure}[ht]
\begin{center}
\include{prismcut}
\end{center}
\caption{Cutting a prism (prismcut)}
\label{prismcut} 
\end{figure}
If we denote a cube by $C$, and a (pentagonal) dumpling by $D$, then the five cuttings can be presented as
$$T+D\quad({\rm triangle}),\quad P+D\quad({\rm quadrilateral}),\quad P+C\quad({\rm quadrilateral}),$$$$ P+P\quad({\rm triangle}),\quad D+D\quad({\rm pentagon}).$$
The last one will appear again in \S\ref{67}.
\subsubsection{A few facts seen from the above cuttings}
Since five planes cut out tetrahedra and prisms, if two chambers cut out by an arrangement of six planes are adjacent along a face, the union is a tetrahedron or a prism. 
Thus though $C$ and $D$ have quadrilateral faces, they are not face to face.

When two $D$'s are face to face along a pentagon, the remaining two pentagons do not share an edge. 
Since there is one pentagon on each plane (see Introduction), and since the group $\mathbb{Z}_6$ acts on the arrangement, we conclude that there are six $D$'s, which are glued to form a solid torus (see \S3.4).  
\subsection{Setting}
If we can choose coordinates on $\mathbb{P}^n$ and the equations of $n+3$ hyperplanes so nicely that the $\mathbb{Z}_{n+3}$-action is clearly seen, it would be nice. 
Though we can not expect this in general, when $n=3$, there is a very good choice (\cite{Mo} and remark after proposition \ref{orbit}).
We work in the real projective space coordinatized by $x:y:z:t$. Our six planes are
$$\begin{array}{lll}
H_x: x=0,& H_y: y=0,& H_z: z=0,\\[3mm]
H^x:h^x=0,& H^y: h^y=0,& H^z: h^z=0,\end{array}$$
where
$$h^x:=y-z-t,\quad h^y:= z-x-t,\quad h^z:= x-y-t.$$
  
Note that if we put $h^t=x+y+z$, then we have 
$$h^y-h^z-h^t=-3x,\quad h^z-h^x-h^t=-3y,\quad h^x-h^y-h^t=-3z.$$
Note also that if we change from $t$ to $-t$, the new arrangement is the mirror image of the original one.

We often work in the Euclidean 3-space coordinatized by
$$(x,y,z)=x:y:z:1,$$
especially when we speak about distance and/or angle, without saying so explicitly.
\subsubsection{A group action and a cubic curve $K$}
\label{cubic}
The six planes admit the transformations 
$$\rho:H_x\to H_y\to H_z\to H_x,\ H^x\to H^y\to H^z\to H^x,$$
$$\sigma: H_x\to H_y, H_y\to H_x, H_z\to H_z,\ H^x\to H^y, H^y\to H^x, H^z\to H^z,\quad\mbox{and}$$
$$\tau:H_a\to H^a\to H_a,\quad (a=x,y,z)$$ 
which are of order 3, 2 and 2, respectively. These are given by the projective transformations
$$
\rho=\left(\begin{array}{cccc}
0&1&0&0\\0&0&1&0\\1&0&0&0\\0&0&0&1
\end{array}\right),\quad 
\sigma=\left(\begin{array}{cccc}
0&\bar1&0&0\\\bar1&0&0&0\\0&0&\bar1&0\\0&0&0&1
\end{array}\right)\quad{\rm and}\quad \tau=\left(\begin{array}{cccc}0&1&\bar1&\bar1\\
\bar1&0&1&\bar1\\1&\bar1&0&\bar1\\1&1&1&0\end{array}\right),$$
where $\bar 1=-1$. Note that $\tau^2=-3I_4$.
In the euclidean space, $\rho$ acts as the $2\pi/3$-rotation around the axis generated by the vector $(1,1,1)$, 
and $\sigma$ the $\pi$-rotation about the axis generated by the vector $(1,-1,0)$. The transformation $\tau$ exchanges the plane 
$$H^{\infty}:x+y+z=0$$ and the plane $H_\infty:t=0$ at infinity. 
They generate the orientation preserving projective transformation group
$$
G \cong\langle \rho,\sigma\rangle\times\langle \tau\rangle\cong 
\mathfrak{S}_3\times\mathbb{Z}_2,
$$
which is a subgroup of $PGL_+(4,\mathbb{Z})$ of order $12$ with center $\langle\tau\rangle=\{1,\tau\}$.  The relations $\sigma^2=1$, $(\rho\tau)^6=1$ and $(\rho\tau\sigma)^2=1$ show that 
$G$ is also isomorphic to the dihedral group $$D_6\cong\langle\rho\tau\rangle\rtimes\langle\sigma\rangle\cong \mathbb{Z}_6\rtimes\mathbb{Z}_2.$$
\begin{prp}
\label{orbit}
There is no invariant plane by $G$. There are two orbits of order two,
$$
\{H^\infty,H_\infty\}\quad\text{and}\quad\{H^{\sqrt{3}},H_{\sqrt{3}}\},
$$
where 
$$
H^{\sqrt{3}}: x+y+z-t\sqrt{3}=0,\quad H_{\sqrt{3}}: x+y+z+t\sqrt{3}=0. 
$$
There is no orbit of order three.
\end{prp}
{\bf proof.} 
If $H$ is invariant by $\tau$, it is different from $H_\infty$. If, in addition, it is invariant by $\rho$ it is orthogonal (in the euclidean space)
to the vector $(1,1,1)$. Such a plane is not invariant by $\sigma$, so that there is no plane invariant by $G$.

Suppose now that $H$ belongs to an orbit of order $2$, say $\{H,H'\}$, different from $\{H^\infty,H_\infty\}$. The group $\langle\rho\rangle\cong\mathbb{Z}_3$
acts necessarily trivially on a set of order $2$ so that $\rho(H)=H$. Then, as above, $H$ is orthogonal to the vector $(1,1,1)$ 
and is not invariant by $\sigma$ so that $\sigma(H)=H'$.
The equation of $H$ writes $x+y+z-\alpha t=0$, and the one of $H'$ writes $x+y+z+\alpha t=0$. Writing that $\tau(H)=H$ 
or $\tau(H)=H'$, we find that $\alpha=\pm\sqrt{3}$.

Now, suppose that $H$ belongs to an orbit of order $3$, then the order two elements $\sigma$ and $\tau$ act trivially on this orbit, so that
$\sigma(H)=\tau(H)=H$, and similarly for $\rho(H)$ and $\rho^2(H)$. In particular $H\neq H_\infty$ and, 
since $\sigma$ is a half-turn rotation about the vector $(1,-1,0)$, either $H$ is orthogonal to the vector $(1,-1,0)$ or $H$ contains 
the line generated by the vector $(1,-1,0)$. 

In the first case, $\rho(H)$ is no longer orthogonal to the vector $(1,-1,0)$ so that $\rho(H)$,
and $\rho^2(H)$ as well, must contain the line generated by the vector $(1,-1,0)$, so that $\rho(H)$, and then $H$ itself, 
is orthogonal to the vector $(1,1,1)$. Therefore $H$ would be invariant, which is impossible.

In the second case, since $H$ is not invariant and then not orthogonal to the vector $(1,1,1)$, the plane $\rho(H)$
doesn't contain the line generated by the vector $(1,-1,0)$, and therefore is orthogonal to the vector $(1,-1,0)$. We
are brought back to the first case.
$\square$

\par\bigskip
Each plane intersects the remaining five; they cut out a pentagon, five triangles and five rectangles. 
See Figures \ref{planeH}, and  \ref{planeH-} for $H_\bullet$ and $H^\circ$, respectively.

The six planes form a Veronese arrangement in the order: 
$$H_x,\quad H^y,\quad H_z,\quad H^x,\quad H_y,\quad H^z,\quad H_x,$$
by which we mean, there is a unique rational cubic curve $K$ osculating these six planes in this order. Note that $\rho$ and $\tau$ respect the order, and $\sigma$ just reverses the order.
Proof goes as follows. 

If the regular hexagon (the projective line coordinatized by $s$) is given by
$$s=0,\quad 1,\quad 3/2,\quad 2,\quad 3,\quad \infty,$$
(here, regular means there is a projective transformation $s\to 3/(3-s)$ sending $0$ to $1$, $1$ to $3/2$, \dots, $\infty$ to $0$)
then the osculating curve is given by
$$K: x=s^3,\quad y=(s-3)^3,\quad z=-8(s-3/2)^3,\quad t=9(s-1)(s-2),$$
which osculates the planes $H_x,\dots,H^z$ at
$$
\begin{array}{lll}K(0)=0:3:\bar3:2, &K(1)=1:\bar8:1:0, &K(\frac32)=\bar3:3:0:2,\\[3mm] K(2)=\bar8:1:1:0,&K(3)=3:0:\bar3:2, &K(\infty)=1:1:\bar8:0.
\qquad \square\end{array}
$$

\noindent{\bf Remark}. The chosen arrangement can be called {\bf selfadjoint} in the following sense. The dual coordinates of the planes are 
$$
H_x=1000,\quad H^y=\bar101\bar1,\quad H_z=0010,\quad H^x=01\bar1\bar1,\quad H_y=0100,\quad H^z=1\bar10\bar1.
$$
For instance the equation of the plane $H^y=\bar101\bar1$ writes $-x+z-t=0$. The six planes turn out to be the points 
$$
P_x=1000,\quad P^y=\bar101\bar1,\quad P_z=0010,\quad P^x=01\bar1\bar1,\quad P_y=0100,\quad P^z=1\bar10\bar1
$$
in the dual space coordinatized by $\check x:\check y:\check z:\check t$. They are on the cubic curve $\check K$ given by
$$
\check x=2(s-1)(s-3/2)(s-2),\quad \check y=-s(s-1),\quad \check z=(s-2)(s-3),\quad \check t=-3(s-1)(s-2),
$$
in this order. Indeed we have
$$
\check K(\infty)=P_x,\quad \check K(0)=P^y,\quad \check K(1)=P_z,\quad \check K(3/2)=P^x,\quad \check K(2)=P_y,\quad \check K(3)=P^z.
$$
The planes supported by three consecutive points
$$
\langle P_x,P^y,P_z\rangle\quad\langle P^y,P_z,P^x\rangle\quad\langle P_z,P^x,P_y\rangle\quad\langle P^x,P_y,P^z\rangle
\quad\langle P_y,P^z,P_x\rangle\quad\langle P^z,P_x,P^y\rangle
$$
happen to be the original ones:
$$
H_x,\quad H^y,\quad H_z,\quad H^x,\quad H_y,\quad H^z.
$$
\par\bigskip
We sometimes code these six planes $H_x,H^y,\dots$ as $H_1,\dots,H_6$. 
Then the transformation of $\mathbb{P}^1$ above induces a projective transformation of $\mathbb{P}^3$ sending $H_j \to H_{j+1}$ mod 6. 
This generates the cyclic group $\langle\rho,\tau\rangle\cong\mathbb{Z}_6$.

\subsubsection{The twenty points}
Convention: $H_{a\cdots}^{b\cdots}:=H_a\cap\cdots\cap H^b\cap\cdots$.
Three planes meet at a point; there are ${6\choose3}=20$ of them.
They are divided into three $G$-orbits:

$$ \text{\ding{108}}\quad P_0=H_{xyz}=0:0:0:1,\quad \bigcirc\quad P^0=H^{xyz}=1:1:1:0;$$$$$$ 
$$\begin{array}{lll}\blacksquare\quad  H^{yz}_x=0:\bar1:1:1,& H^{zx}_y=1:0:\bar1:1,& H^{xy}_z=\bar1:1:0:1,\\[3mm]
\square\quad H_{yz}^x=1:0:0:0,& H_{zx}^y=0:1:0:0,& H_{xy}^z=0:0:1:0;\\[7mm]
\bullet\quad H_{xy}^x=0:0:\bar1:1,& H_{yz}^y=\bar1:0:0:1,& H_{zx}^z=0:\bar1:0:1,\\[3mm] 
\circ\quad H^{xy}_x=0:2:1:1,& H^{yz}_y=1:0:2:1,& H^{zx}_z=2:1:0:1;\\[3mm]
\bullet\quad H_{yx}^y=0:0:1:1,& H_{zy}^z=1:0:0:1,& H_{xz}^x=0:1:0:1,\\[3mm]
\circ\quad H^{yx}_y=\bar2:0:\bar1:1,& H^{zy}_z=\bar1:\bar2:0:1,& H^{xz}_x=0:\bar1:\bar2:1,\end{array}$$
$$$$ 
where $\bar1=-1, \bar2=-2$. 
Note that the distances from the vertices to the origin $P_0$ are either 
$$\text{\ding{108}}:0,\qquad \bullet:1,\quad\blacksquare:\sqrt2,\quad\circ:\sqrt5,\qquad \square,\ \bigcirc:\infty.$$
In the $1\cdots6$ -coding, the three orbits above are represented by
$$135,\quad 123, \quad 124,$$
respectively, where $135$ stands for $H_1\cap H_3\cap H_5$.

\subsubsection{An invariant quadratic form}
While there are many $G$-invariant quadratic forms in $(x,y,z,t)$, if we ask the zero set passes through the six vertices marked $\square$ and $\blacksquare$, then it is (up to scalar multiplication) given by
$$Q=xy+yz+zx+t^2.$$
The surface $Q=0$ and the curve $K$ do not meet: in fact, substituting the expression of the curve in $s$ into $Q$, we have $-15(s^2-3s+3)^3$.

By identifying the space $\mathbb P^3$ and its dual, we regard the cubic curve $\check K$ live in our space, that is, the curve is defined by
$$x=2(s-1)(s-3/2)(s-2),\quad  y=-s(s-1),\cdots$$
Then this is the unique rational cubic curve passing through the six points marked $\square$ and $\blacksquare$:
$$\check K(\infty)=H^x_{yz},\quad \check K(0)=H^{zx}_y,\quad \check K(1)=H^z_{xy},\quad \check K(3/2)=H^{yz}_x,\quad \check K(2)=H^y_{zx},\quad \check K(3)=H^{xy}_z.$$
This curve $\check K$ is on the surface $Q=0$.  

\subsubsection{The twenty six chambers}
The six planes in $\mathbb{P}^3$ cut out 
\par\medskip
\begin{center}two cubes \ $C$,\ twelve prisms \  $P$,\
six tetrahedra \ $T$,\ six dumplings \  $D$.
\end{center}
The cubic curve $K$ stays in the dumplings, and osculates each pentagon. 
The quadratic surface $Q=0$ lies in the union of the tetrahedra and the dumplings (see Figures \ref{planeH}, \ref{planeH-}).

\subsubsection{ The six planes}
For the planes $H^x,H^y,H^z$, the intersection with other planes are shown in  Figure \ref{planeH-}.
\begin{figure}[h] 
\begin{center}
\include{planeH-}
\end{center}
\caption{ The planes $H^\circ$ (planeH-)}
\label{planeH-}
\end{figure}
\noindent
Here $TP, DP,\dots$ stand for $T\cap P, D\cap P,\dots.$ If the plane $H^x$, defined by $z=y-1$, is coordinatized by $(y,x)$, 
then the dotted curve is the conic (a hyperbola) 
$$
Q=xy+y(y-1)+(y-1)x+1=\left(y-\frac12\right)\left(y+2x-\frac12\right)+\frac34=0.
$$

For the planes $H_x,H_y,H_z$, the intersection with other planes are shown in  Figure \ref{planeH}.
\begin{figure}[h]
\begin{center}
\include{planeH}
\end{center}
\caption{ The planes $H_\bullet$ (planeH)}
\label{planeH} 
\end{figure}
\noindent
If the plane $H_x$, defined by $x=0$, is coordinatized by $(y,z)$, then the dotted curve is the conic (hyperbola) $Q=yz+1=0$, and the point $\heartsuit$ is the point  $K(0)$ where the curve $K$ kisses the plane $H_x$. 

\subsubsection{Bottom and top}
It is convenient to consider the triangle with the three vertices $\blacksquare$ in the plane $H^\infty: x+y+z=0$ as the ground floor, and the triangle with the three vertices $\square$ in the plane $H_\infty$ at infinity as the ceiling. The planes $H^x,H^y,H^z$ are orthogonal to the ground. The intersection (lines) of the six planes with the floor and the ceiling are shown in Figure \ref{bt}. The dotted curves are the intersections with the surface $Q=0$.
\begin{figure}[h]
\begin{center} 
\include{bottom-top} 
\end{center}
\caption{ Bottom and top (bottom-top)}
\label{bt}
\end{figure}
%
%
\subsubsection{Octant chambers}\label{subsectionOctant}
The finite space coordinatized by
$$(x,y,z)=x:y:z:1$$
is divided by the three planes $H_x,H_y$ and $H_z$ into eight chambers.
We denote them as 
$$(+++)=\{x\ge0,y\ge0,z\ge0\},\quad (++-)=\{x\ge0,y\ge0,z\le0\},\dots.$$
In this section, we see how these chambers are cut by the three planes 
$H^x,H^y$ and $H^z$. In the chambers $(+++)$ and $(---)$, the happenings are similar, and in the remaining six chambers, similar things happen. Set
$$(+++;---)=\{(x,y,z)\in(+++)\mid h^x\le0,h^y\le0,h^z\le0\},...$$
If we write these as $(\e_1\e_2\e_3;\eta_1\eta_2\eta_3)$, they permit the $\mathbb{Z}_3$-action $$(\e_1\e_2\e_3;\eta_1\eta_2\eta_3)\to(\e_2\e_3\e_1;\eta_2\eta_3\eta_1).$$
\par\medskip\noindent{\bf In the chamber $(+++)$:}
There are four $\mathbb{Z}_3$-orbits represented by
$$\begin{array}{ll}(+++;+++)=\emptyset,\quad &(\dot+\ \dot+\ \dot+;\dot-\ \dot-\ \dot-):{\ \rm cube}\\[3mm]
(\dot++\dot+;\dot+\ \dot-\ \dot-): {\ \rm pri\ }P',\quad
&(+\ \dot+\ +;-\ \dot+\ \dot+): {\ \rm sm\ }P''.\end{array}$$
Here effective ones are marked by dots.
For example, the last one is defined by $y\ge0,h^y\ge0$ and $h^z\ge0$; these three inequalities imply the other ones $x\ge0,z\ge0$ and $h^x\le0$. 
The number of dots corresponds to the number of walls.

As a whole, there are a cube and three $P'$ and three $P''$ 
(see Figure \ref{octants}). 
One $P'$ and another $P''$ in the opposite chamber $(---)$ are glued along 
the plane at infinity forming a (full) prism (see Figure \ref{inftyprism80}).
\begin{figure}[h] 
\begin{center}  
\include{inftyprism80}   
\end{center}
\caption{Unbounded prism $P=P'\cup P''$ (inftyprism80)}
\label{inftyprism80} 
\end{figure}

If we cut this chamber by a big sphere centered at the origin, the intersection with the chamber is surrounded by three arcs ($x=0, y=0, z=0$), 
and the arc-triangle is cut by the three lines ($h^x=0,h^y=0,h^z=0$). Around a triangle (cube), there are a triangle ($P''$), a pentagon ($P'$), a triangle ($P''$),... (see Figure \ref{sphere1} and Figure \ref{sphere3} left).
\par\medskip\noindent{\bf In the chamber $(---)$:}
There are four $\mathbb{Z}_3$-orbits represented by
$$\begin{array}{ll}(---;+++)=\emptyset,\quad &(\dot-\dot-\dot-;\dot-\dot-\dot-):{\ \rm cube}\\[3mm]
(\dot-\dot--;\dot+\dot-\dot-): {\ \rm pri\ }P',\quad
&(--\dot-;-\dot+\dot+): {\ \rm sm\ }P''.
\end{array}$$ 

\begin{figure}
  \begin{center} 
  \includegraphics[width=110mm]{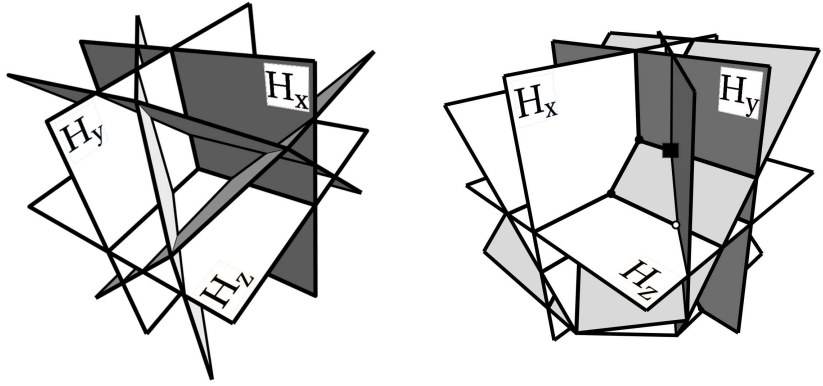}
  \end{center} 
  \caption{Chambers $(+++)$ and $(++-)$ (octants)} 
  \label{octants}
\end{figure}

\par\medskip\noindent{\bf In the chamber $(++-)$:}
This chamber does not admit the group action. Since $x\ge0$ and $z\le0$ imply $z-x\le0$, and so $h^y=z-x-1<0$, we have $(++-;*+*)=\emptyset.$ 
Among the remaining four, there is a unique compact one:
$$(\dot+\dot+\dot-;\dot--\dot-): \ {\rm prism}.$$
(See Figures \ref{octants}(right) and \ref{octant70}(left).) The others are
$$(+\dot+\dot-;\dot--\dot+): \ {\rm tetrahedron},$$
and 
$$(\dot+\dot+\dot-;\dot+-\dot-): \ {\rm dump}\ D',\quad(+\dot+\dot-;\dot+-\dot+): \ {\rm ling}\ D''.$$
One $D'$ (resp. $D''$) and another $D''$ (resp. $D'$) in the opposite chamber $(--+)$ are glued along the plane at infinity forming a (full) dumpling (see Figure \ref{octant70}).  

If we cut this chamber by a big sphere centered at the origin, the chamber is surrounded by three arcs ($x=0, y=0, z=0$), and the arc-triangle is cut by the two lines ($h^x=0,h^z=0$). There are a triangle (tetrahedron) and two quadrilaterals, which are sections of $D'$ and $D''$ (see Figure \ref{sphere1} and Figure \ref{sphere3} right).
\par\medskip\noindent{\bf In the chamber $(--+)$:}
Since $y\le0$ and $z\ge0$ imply $y-z\le0$, and so $h^x=y-z-1<0$, we have
$(--+;+**)=\emptyset$.  
Among the remaining four, there is a unique compact one:
$$(\dot-\dot-\dot+;-\dot-\dot-): \ {\rm prism}.$$
The others are
$$(\dot--\dot+;-\dot-\dot+): \ {\rm tetrahedron},$$
and 
$$(\dot-\dot-\dot+;-\dot+\dot-): \ {\rm dump}\ D',\quad(\dot--\dot+;-\dot+\dot+): \ {\rm ling}\ D''.$$

\begin{figure}[h] 
\begin{center}  
\include{octant70} 
\end{center}
\caption{(right) A dumpling $D$ is cut by the plane at infinity into two parts dump $D'$ and ling $D''$. (left) They are in octants of type $(++-)$. The faces of $D$ have names
$1' + 1'', 2, 3'+3'', 4, 5, 6'+6''$.
In the octant, a wall is $6'$ seen from $D'$, and is $6''$ seen from $D''$.  (octant70)}
\label{octant70} 
\end{figure}

\newpage
\subsection{The solid torus made of the two cubes and the twelve prisms}
The two cubes kiss at the two vertices $P_0$ and $P^0$, which are opposite 
vertices of each cube. Around the two cubes are twelve prisms, 
forming with the two cubes a solid torus. We explain how they are situated.

\subsubsection{The two cubes and the six prisms around $P_0$}

\begin{figure}[h] 
  \begin{center}
  \includegraphics[width=140mm]{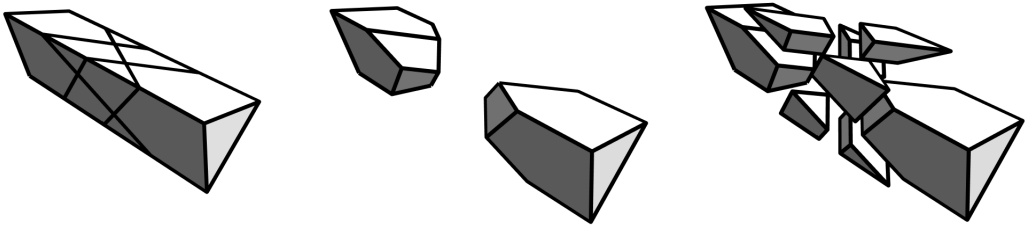}
  \end{center} 
  \caption{(left) The big-cylinder. (middle) The two truncated cubes slightly moved away from each other. (right) 
  Exploded view of the two cubes and the six prisms.(amy-fig13.eps)} 
  \label{infiniteprism}
\end{figure}

The two cubes and the six prisms 
around $P_0$, form in the finite (Euclidean) space an infinitely long triangular cylinder called {\bf the big-cylinder} 
and bounded by $H^x,H^y,H^z$. It is defined 
by the inequations 
(figure \ref{infiniteprism}):
$$
h^x\leqslant0,\quad h^y\leqslant0,\quad h^z\leqslant0. 
$$

\begin{figure}[h]
  \begin{center}
  \includegraphics[width=140mm]{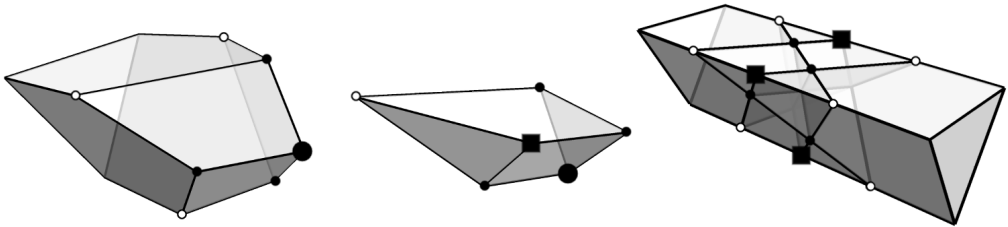}
  \end{center} 
  \caption{The big-cylinder made by two truncated cubes and six compact prisms (right). A truncated cube (left) and a prism (center) 
are shown with their vertices.} 
  \label{truncatedcube}
\end{figure}

\noindent Part of the pictures in Figure \ref{infiniteprism} is enlarged in 
Figure \ref{truncatedcube}.

\subsubsection{The small cylinder}

\begin{figure}
  \begin{center}
  \includegraphics[width=120mm]{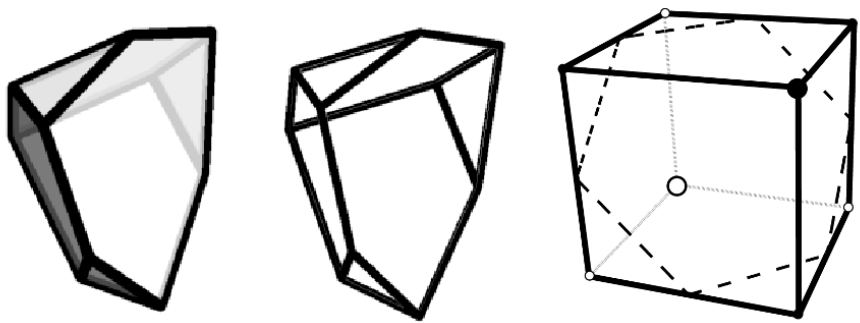}
  \end{center} 
  \caption{The heptahedron $C_+$ (left) and its 1-skeleton (center). Cutting a cube into two heptahedra (right)} 
  \label{heptahedron}
\end{figure}

Considering the planes $H^{\sqrt{3}}$ and $H_{\sqrt{3}}$ defined in proposition \ref{orbit},
we have $\tau(H_{\sqrt{3}})=H^{\sqrt{3}}$ where $\tau$ is the involution
defined in \ref{cubic}. 
Let us denote by $C_+$ the heptahedron shown in figure \ref{heptahedron} and defined by 
$$
\frac{x}{t}\geqslant0,\quad \frac{y}{t}\geqslant0,\quad \frac{z}{t}\geqslant0,\quad \frac{x}{t}-\frac{y}{t}\leqslant1,
\quad\frac{y}{t}-\frac{z}{t}\leqslant1,\quad\frac{z}{t}-\frac{x}{t}\leqslant1,
\quad\frac{x+y+z}{t}\leqslant\sqrt{3}.
$$
 
\begin{figure}
  \begin{center}
  \includegraphics[width=120mm]{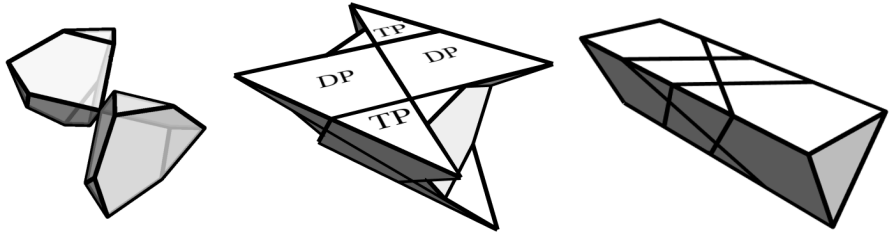}
  \end{center} 
  \caption{Heptahedrons $C_+$ and $C_-$ kissing at $P_0$ (left). 
  The small cylinder $SC$ and two quadrilaterals $DP$ and two triangles $TP$ on its boundary (middle). 
  The triangular prism homeomorphic to $SC$ (right).} 
  \label{smallcylinder}
\end{figure}

\noindent We set $C_-:=\sigma(C_+)$. The faces of the heptahedron $C_+$ are three triangles, 
three pentagons and one hexagon. 
We consider the solid $SC$ (called {\bf the small cylinder}) given by the union of the two 
heptahedrons and the six prisms. It is homeomorphic to a triangular prism by a homeomorphism 
preserving the generating lines (see figure \ref{smallcylinder}).

\subsubsection{The prism torus}
The $G$-orbit of one of the six compact prisms around $P_0$ is given by twelve prisms,
the six extra ones being located around $P^0$. Now the image of the small cylinder $SC$
by $\sigma$, or equivalently, the $G$-orbit of one prism and the heptahedron $C_+$, 
form the union of the big cylinder and the six prisms around $P^0$. 
These two solid cylinders $SC$ and $\sigma(SC)$ are attached by $\sigma$ along their 
bases and therefore form a solid torus. Let us call it the {\bf prism torus}.
Thus, this prism torus is decomposed into sixteenteen cells, four heptahedrons and twelve prisms,
on which acts the group $G$.
 
\subsubsection{The boundary of the prism torus}\label{dprismT}
The prism torus consists also of the two cubes and the twelve prisms. 
Note that this torus is a tubular neighborhood of the projective line joining 
$P_0$ and $P^0$. Therefore, as mentioned in the introduction of section \ref{6planes}, the 
complementary of this torus is a solid torus as well.

\noindent
The above construction shows that the boundary of the torus consists of the faces 
of the prisms on the lateral boundary on the small cylinder $SC$ and their images 
by $\sigma$. In particular, no face of the cubes belongs to the boundary 
of the prism torus. 

\noindent The boundary of the prism torus on each plane consists of two quadrilaterals $DP$ and two triangles $TP$ (see figure \ref{smallcylinder}). 
 They are shown more explicitly the pattern below:
\scriptsize
$$
\begin{array}{ccc}
\begin{array}{ccccccc}
\bullet&     -    &\square&-&\bullet&&\\
       &\backslash&   |   & & |     &&\\
       &          &\circ  &-&\circ  &&\\
       &          &   |   & & |     &\backslash&\\
       &          &\bullet&-&\square&   -      &\bullet
\end{array}
&\qquad&
\begin{array}{ccccccc}
\circ&     -    &\blacksquare&-&\circ&&\\
       &\backslash&   |   & & |     &&\\
       &          &\bullet&-&\bullet&&\\
       &          &   |   & & |     &\backslash&\\
       &          &\circ&-&\blacksquare&   -      &\circ
\end{array}
\end{array}
$$\normalsize

\begin{figure}[ht]
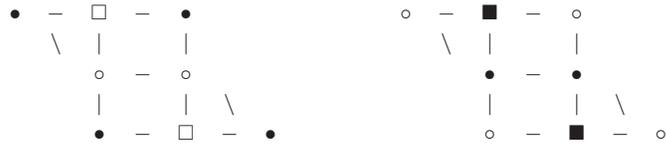

\begin{center}
\include{planeH2-70}
\end{center}
\caption{ The boundary of the prism torus consists of the faces of the prisms (planeH2-70)}
\label{planeH2-70}
\end{figure}

The twelve quadrilaterals and the twelve triangles tessellate the boundary torus as below. To indicate the identification, four copies of the plane number 1 are shown in addition. After identification, there are 
$$3\  \square,\quad 3\ \blacksquare,\quad 6\ \bullet,\quad 6\ \circ;$$
note that among the twenty vertices, only $P^0$ and $P_0$ are missing.
\par\medskip
\scriptsize
$$\begin{array}{cccccccccccccccccc}
   &   &   &   &   &   &   &   &   &   &   &\bu& - &\sq& - &\bu&   &\\
   &   &   &   &   &   &   &   &   &   &   &   &\ba& | &   & | &   &\\
   &   &   &   &   &   &   &   &   &\ci& - &\bl& - &\ci& 1 &\ci&   &\\
   &   &   &   &   &   &   &   &   &   &\ba& | &   & | &   & | &\ba&\\
   &\ci& - &\sq& - &\bu&   &\bu& - &\sq& - &\bu& 2 &\bu& - &\sq& - &\bu\\
   &   &\ba& | &   & | &   &   &\ba& | &   & | &   & | &\ba&   &   &\\
   &   &   &\ci& 1 &\ci& - &\bl& - &\ci& 3 &\ci& - &\bl& - &\ci&   &\\
   &   &   & | &   & | &\ba& | &   & | &   & | &\ba&   &   &   &   &\\
   &   &   &\bu& - &\sq& - &\bu& 4 &\bu& - &\sq& - &\bu&   &   &   &\\
   &   &   &   &\ba& | &   & | &   & | &\ba& | &   & | &   &   &   &\\
   &\ci& - &\bl& - &\ci& 5 &\ci& - &\bl& - &\ci& 1 &\ci&   &   &   &\\
   &   &\ba& | &   & | &   & | &\ba&   &   & | &   & | &\ba&   &   &\\
 - &\sq& - &\bu& 6 &\bu& - &\sq& - &\bu&   &\bu& - &\sq& - &\bu&   &\\
\ba& | &   & | &   & | &\ba&   &   &   &   &   &   &   &   &   &   &\\
   &\ci& 1 &\ci& - &\bl& - &\ci&   &   &   &   &   &   &   &   &   &\\
   & | &   & | &\ba&   &   &   &   &   &   &   &   &   &   &   &   &\\
   &\bu& - &\sq& - &\bu&   &   &   &   &   &   &   &   &   &   &   &\\
\end{array}$$
\normalsize

%
\subsection{The solid torus made of the six dumplings}
Each of the six dumplings is glued along the two pentagonal faces with two others. The six dumplings form a solid torus; let us call it the {\bf dumpling torus}. This section is devoted to understand it.
\subsubsection{A dumpling}
A dumpling has two rectangular faces, two triangular faces, and two pentagonal faces (see Figure \ref{planeH3-70}),  which share an edge with two vertices marked white and black squares. This edge will be called a {\bf special edge} (see Figure \ref{gyoza170}). 

\begin{figure}[ht] 
\begin{center}
\include{planeH3-70}
\end{center}
\caption{ Two kinds of pentagons bounding the dumplings (planeH3-70)}
\label{planeH3-70}
\end{figure}
\begin{figure}[ht]
\begin{center}
\include{gyoza170}
\end{center}
\caption{ A dumpling, and two dumplings glued along a pentagon (gyoza170)}
\label{gyoza170}
\end{figure}

Each dumpling $D$ is cut into two by the plane at infinity: the part with two vertices $\bullet$ is denoted by $D'$, and the other by $D''$. They are shown in Figure \ref{octant70}.

%

\subsubsection{The dumpling torus}
The six dumplings glued along their pentagonal faces form the dumpling torus.
Figure \ref{pentagons70} shows the six pentagons glued along the special edges forming a circle.
In the figure, special edges are shown by thick segments.
\begin{figure}[ht]
\begin{center}
\include{pentagons70} 
\end{center}
\caption{ Theix pentagons glued along the special edges (pentagons70)}
\label{pentagons70}
\end{figure} 

\begin{figure}[ht]
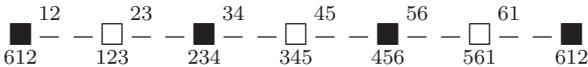
 
\begin{center}
\include{pileofG1} 
\end{center}
\caption{Piles of six dumplings glued along their pentagonal faces (pileofG1)}
\label{pileofG1}
\end{figure}

Remember that the planes
$$H_x, H^y, H_z, H^x, H_y, H^z, H_x\quad\mbox{\rm are\ numbered}\quad 1,\dots,6;$$ and the curve $K$ lives inside the dumpling torus, and touches six pentagons in this order.
In Figure \ref{pileofG1}, the six dumplings are glued along their pentagonal faces. The smallest pentagon on the top and the biggest one at the bottom are on the plane number 6; they should be identified after $2\pi/5$ turn. The pentagon on the plane 2 is next to the top, and so on. The edges are labeled as 12, 23, ...; for example, 12  indicates the intersection of two pentagons on the planes 1 and 2.  The edges labeled by two consecutive numbers are the special edges.  The special edges form a circle:
$$\underset{612}{\blacksquare}\overset{12}{--}\underset{123}{\square}\overset{23}{--}\underset{234}{\blacksquare}\overset{34}{--}\underset{345}{\square}\overset{45}{--}\underset{456}{\blacksquare}\overset{56}{--}\underset{561}{\square}\overset{61}{--}\underset{612}{\blacksquare}$$
\newpage  
One can see in Figure \ref{pileofG1}, in front, the quadrilateral consisting of two triangles and two rectangles with edges labeled 12 at the top and 61 in the bottom, and segments 31, 41 and 51 are inside; this rectangle together with the pentagon $12-13-14-15-61$ is on the plane 1.
\subsubsection{Some curves on the dumpling torus}
Trace the edge 51 on this rectangle along non-special edges, and we have the close curve:
$$\underset{561}{\square}\overset{51}{--}\bullet\overset{52}{--}\bullet\overset{53}{--}\underset{345}{\square}\overset{35}{--}\bullet\overset{36}{--}\bullet\overset{31}{--}\underset{123}{\square}\overset{13}{--}\bullet\overset{14}{--}\bullet\overset{15}{--}\underset{561}\square.$$
Trace the edge 41 on this rectangle along non-special edges, and we have the close curve: 
$$\underset{612}{\blacksquare}\overset{62}{--}\circ\overset{63}{--}\circ\overset{64}{--}\underset{456}{\blacksquare}\overset{46}{--}\circ\overset{41}{--}\circ\overset{42}{--}\underset{234}{\blacksquare}\overset{24}{--}\circ\overset{25}{--}\circ\overset{26}{--}\underset{612}\blacksquare.$$
The other edges form a single closed curve:
$$\underset{561}{\square}\overset{61}{--}\underset{4}\circ\overset{61}{--}\underset{3}\bullet\overset{61}{--}\underset{612}{\blacksquare}\overset{12}{--}\underset{5}\bullet\overset{12}{--}\underset{4}\circ\overset{12}{--}\underset{123}{\square}\overset{23}{--}\underset{6}\circ\overset{23}{--}\underset{5}\bullet\overset{23}{--}\underset{234}\blacksquare$$
$$\underset{234}{\blacksquare}\overset{34}{--}\underset{1}\bullet\overset{34}{--}\underset{6}\circ\overset{34}{--}\underset{345}{\square}\overset{45}{--}\underset{2}\circ\overset{45}{--}\underset{3}\bullet\overset{45}{--}\underset{456}{\blacksquare}\overset{56}{--}\underset{2}\bullet\overset{56}{--}\underset{3}\circ\overset{56}{--}\underset{561}\square.$$
Boundary of a(ny) pentagon; they will be called $m_D$ later:
$$\blacksquare--\circ--\circ--\blacksquare--\square--\blacksquare,\quad \square--\bullet--\bullet--\square--\blacksquare--\square.$$
The diagonal of the front quadrilateral in Figure \ref{pileofG1} connecting $\underset{612}\blacksquare$ crossing the edges $31,41$ and $51$; this will be called $L$ later.
\subsubsection{Another description of the dumpling torus}
The dumpling torus presented as the pile of dumplings above can be understood also as follows, which may help the understanding:  Consider the cylinder made by six copies of a pentagon times the interval $[0,1]$, glued along pentagonal faces; the top and the bottom with $2\pi/5$ twist. The side of the cylinder is tessellated by  $5\times 6$ rectangular tiles. Now choose six squares diagonally arranged (thanks to the twist, one can make it consistently) and compress these squares vertically, then you end up with the pile of dumplings above. In Figure \ref{chochin}, it is shown in particular, how the dumpling with pentagonal faces 5 and 6 is made by a pentagon times the interval.
\begin{figure}[ht]
\begin{center} 
\include{chochin}
\end{center}
\caption{Compressing the pile of six pentagonal prisms to the dumpling torus (chochin)}
\label{chochin}
\end{figure}


\subsubsection{The boundary of the dumpling torus and the line $L$}
The boundary of the dumpling torus on each plane is shown in Figures \ref{planeH} and \ref{planeH-}, and more explicitly again in Figure \ref{planeH4-70}. It consists of two quadrilaterals $DP$ and two triangles $DT$. This is the front rectangle shown in Figure \ref{pileofG1}, and a column of Figure \ref{chochin} left. We show it again below by a diagram which will fit the diagram of the boundary of the prism torus shown in \S \ref{dprismT}. 
\scriptsize%
$$
\begin{array}{ccc}
\begin{array}{cccc}
 &   &   &\bl\\
 &   & / & | \\
 &\sq& - &\bu\\
 & | &   & | \\
 &\ci& - &\ci\\
 & | &   & | \\
 &\bu& - &\sq\\
 & | & / &   \\
 &\bl&   &   
\end{array}
&\qquad\qquad\qquad\qquad&
\begin{array}{cccc}
 &   &   &\sq\\
 &   & / & | \\
 &\bl& - &\ci\\
 & | &   & | \\
 &\bu& - &\bu\\
 & | &   & | \\
 &\ci& - &\bl\\
 & | & / &   \\
 &\sq&   &   
\end{array}
\end{array}
$$\normalsize

Note that in the left diagram, the two opposite $\bl$ stands for the same vertex; in the right diagram, two opposite $\sq$ stands for the same vertex. The diagonal joining two opposite $\bl$ in the left (in Figure \ref{planeH4-70} left, the line $x-y+2=0$) and the diagonal joining two opposite $\sq$ in the right (in Figure \ref{planeH4-70} right, a vertical line $y=1/2$) are closed curves, which are homotopically equivalent on the boundary of the dumpling torus; Let us call one of them $L$.

\normalsize
\begin{figure}[ht]
\begin{center}
\include{planeH4-70}
\end{center}
\caption{ Boundary of the dumpling torus on each plane (planeH4-70)}
\label{planeH4-70}
\end{figure}

We now show the six such forming the boundary of the dumpling torus. The twelve quadrilaterals and the twelve triangles tessellate the boundary torus as below. To indicate the identification, four copies of the plane number 1 are shown. After identification, there are 
$$3\  \square,\quad 3\ \blacksquare,\quad 6\ \bullet,\quad 6\ \circ;$$
note that among the twenty vertices, only $P^0$ and $P_0$ are missing.
\scriptsize
$$\begin{array}{cccccccccccccccccc}
   &   &   &   &   &   &   &   &   &   &   &   &   &   &   &\bl&   &\\     
   &   &   &   &   &   &   &   &   &   &   &   &   &   & / & | &   &\\     
   &   &   &   &   &   &   &   &   &   &   &   &   &\sq& - &\bu&   &\\
   &   &   &   &   &   &   &   &   &   &   &   & / & | &   & | &   &\\
   &   &   &   &   &\bl&   &   &   &   &   &\bl& - &\ci& 1 &\ci&   &\\
   &   &   &   & / & | &   &   &   &   & / & | &   & | &   & | &   &\\
   &   &   &\sq& - &\bu&   &   &   &\sq& - &\bu& 2 &\bu& - &\sq&   &\\
   &   &   & | &   & | &   &   & / & | &   & | &   & | & / &   &   &\\
   &   &   &\ci& 1 &\ci&   &\bl& - &\ci& 3 &\ci& - &\bl& - &\ci&   &\\
   &   &   & | &   & | & / & | &   & | &   & | & / & | &   &   &   &\\
   &   &   &\bu& - &\sq& - &\bu& 4 &\bu& - &\sq& - &\bu&   &   &   &\\
   &   &   & | & / & | &   & | &   & | & / & | &   & | &   &   &   &\\
   &   &   &\bl& - &\ci& 5 &\ci& - &\bl& - &\ci& 1 &\ci&   &   &   &\\
   &   & / & | &   & | &   & | & / &   &   & | &   & | &   &   &   &\\
   &\sq& - &\bu& 6 &\bu& - &\sq& - &\bu&   &\bu& - &\sq&   &   &   &\\
   & | &   & | &   & | & / &   &   &   &   & | & / &   &   &   &   &\\
   &\ci& 1 &\ci& - &\bl& - &\ci&   &   &   &\bl&   &   &   &   &   &\\
   & | &   & | & / &   &   &   &   &   &   &   &   &   &   &   &   &\\
   &\bu& - &\sq&   &   &   &   &   &   &   &   &   &   &   &   &   &\\
   & | & / &   &   &   &   &   &   &   &   &   &   &   &   &   &   &\\
   &\bl&   &   &   &   &   &   &   &   &   &   &   &   &   &   &   &
\end{array}$$
\normalsize
\subsection{Tetrahedra}
Let us compare the (tessellated) boundary of the prism torus and that of the dumpling torus.
They share the quadrilateral faces. 
The two solid tori almost fill the space; six tetrahedra are the notches. Each tetrahedron has two triangular faces adjacent to prisms 
and other two triangular faces adjacent to dumplings:
\scriptsize
\par\medskip
$$\begin{array}{ccccccccccccccccc}
   &\ci& - &\bl&   & \qquad  \qquad &\ci& - &\bl&   & \qquad  \qquad  &   &\ci& - &\bl&\\
  \partial & | &\times& | &   &      =          & | &\ba& | &   &      \cup           &   & | & {\bf /} & | &\\
   &\sq& - &\bu&   &               &\sq& - &\bu&   &                 &   &\sq& - &\bu&
\end{array}$$\normalsize

\par\medskip
\par
\begin{figure}
  \begin{center}
  \includegraphics[width=120mm]{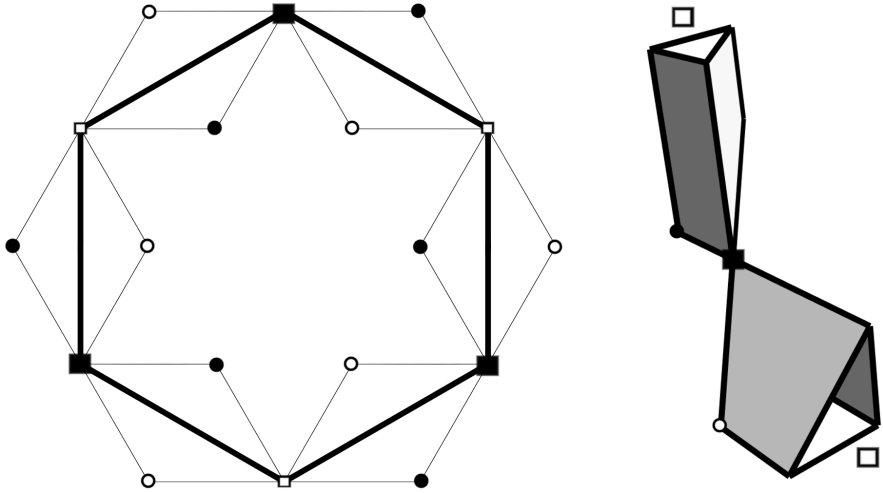}
  \end{center} 
  \caption{The rosary of tetrahedra and two tetrahedra extracted from  Figure \ref{octants}(right)}
  \label{tetra1}
\end{figure}

In Figures \ref{octants}(right) and \ref{octant70}(left), we see that each tetrahedron has bounded base 
(triangle $(\bullet,\circ,\blacksquare)\subset H^\circ$) on the big-cylinder, an infinitely long face on the prism torus, 
and two infinitely long faces on the dumpling torus; intersection of the last two faces is a special edge. 
The six tetrahedra form a rosary (see Figure \ref{tetra1}(left) for a combinatorial idea).

The curve $\check K$, introduced in \S 3.2.3 lives in the union of the six tetrahedra.

\subsection{Intersection with a big sphere}

The aim of this section is to visualize what happens near the plane at infinity $H_{\infty}$ by looking at the intersection with a big sphere.

\noindent The sphere at infinity is shown in Figure \ref{sphere0}; the plane at infinity (see Figure \ref{bt}) is obtained by identifying antipodal points. 
The three circles represent (the intersections with the planes)
$H_x,H_y,H_z,$ showing an octahedron, and the three lines represent 
(the intersections with the planes) $H^x,H^y,H^z$. They cut out 24 triangles. Each triangle is a section of a prism or a dumpling. 
Dotted circles represent (the intersection with) the surface $Q=0$.
\begin{figure}[ht]
\begin{center}
\include{sphere0}
\end{center}
\caption{ Intersection of the chambers with the sphere at infinity (sphere0)}
\label{sphere0}
\end{figure}

A very big sphere is shown in Figure \ref{sphere1}.  The three circles (showing an octahedron) and the three lines represent as in the previous figure. The triangle appeared at the center is the section of a cube; the other cube situates out of the picture. The six small triangles are the sections of the tetrahedra.

\begin{figure}[ht] 
\begin{center} 
\include{sphere1}
\end{center} 
\caption{ Intersection of the chambers with a big sphere (sphere1)}
\label{sphere1}
\end{figure}
Each prism $P$ is the union of two parts $P'$ and $P''$ as is shown in Figure \ref{inftyprism80}, glued along the shaded triangle. Now the shaded triangle is thickened; the sections of $P'$ and $P''$ with the big sphere are a triangle and a pentagon, respectively. (See Figure \ref{sphere3} left.)
\par\medskip
Each dumpling $D$ is the union of two parts $D'$ and $D''$ as is shown in Figure \ref{octant70}, glued along the shaded triangle. Now the shaded triangle is thickened; the sections of $D'$ and $D''$ with the big sphere are both quadrangle. (See Figure \ref{sphere3} right.)

\begin{figure}[ht]
\begin{center}
\include{sphere3}
\end{center}
\caption{ Enlargement of two faces of the octahedron (sphere3)}
\label{sphere3}
\end{figure}

Intersection with the big sphere with the surface $Q=0$ is shown in Figure \ref{sphere2}; 
which shows that the surface does not path through any prisms nor cubes, 
and that the surface can be considered as an approximation of (the boundary of) the prism torus.
\begin{figure}[ht]
\begin{center}
\include{sphere2}
\end{center}
\caption{ Enlargement around the intersection with a tetrahedron together with that of the surface $Q=0$ (sphere2)}
\label{sphere2}
\end{figure}
\clearpage
\section{Seven hyperplanes in the 4-space}
Seven hyperplanes in general position in $\mathbb{P}^4$ cut out a unique chamber, bounded by the seven hyperplanes, 
stable under the action of the cyclic group $\mathbb{Z}_7$;  let us call it the central chamber CC. 
This action can be projective if the hyperplanes are well-arranged. The intersection of CC and a hyperplane is a dumpling. 

So the boundary of CC is a 3-sphere $\mathbb{S}^3$ tessellated by seven dumplings. 
We would like to know how they are arranged, especially the $\mathbb{Z}_7$-action on the tessellation. 
In \cite{CY}, the seven dumplings in $\mathbb{S}^3$ is shown (see Figure \ref{7facebody70}), in which the $\mathbb{Z}_7$-action can be hardly seen.
\par\smallskip
We first observe the tessellation. We label the dumplings as $D_1,D_2,\dots,D_7$ so that the dumplings $D_k$ and $D_{k+1}$, modulo 7, 
share a pentagon, and the group action induces a transformation $D_k\to D_{k+1}$ modulo 7. The special edges form a circle $C$. 
(Recall that a dumpling is bounded by two pentagons, two triangles and two quadrilaterals, 
and that the intersection of the two pentagons is called the special edge.) If we remove seven pentagons $D_k\cap D_{k+1}$ 
from the 2-skeleton of the tessellation, it remains a M\"obius strip $M$ with the curve $C$ as the boundary. 
\par\noindent
We find, in the 2-skeleton, a subskeleton homeomorphic to a disc which bound the curve $C$. 
This implies that $C$ is unknotted, and that $M$ is twisted only by $\pm\pi.$
\par\smallskip
We give another description of the tessellation so that the action of $\mathbb{Z}_7$ can be seen; 
we start from a vertical pentagonal prism, horizontally sliced into seven thinner prisms. 
We collapse vertically the seven rectangular faces diagonally situated on the boundary of the prism; 
this process changes each thinner prism into a dumpling. Identifying the top and the bottom after a suitable rotation, 
we get a solid torus $ST$ made of seven dumplings, so that the solid torus $ST$ admits a $Z_7$-action. 
Consider a usual solid torus $UT$ in our space $\mathbb{S}^3$ and put $ST$ outside of $UT$; the boundary of $ST$, which is also the boundary of $UT$, is a torus $T$. 
We collapse $UT$ by folding the torus $T$ into a M\"obius strip, which can be identified with $M$, whose boundary can be identified with $C$.
We will see that the result, say $X_7$, of this collapsing is homeomorphic to $\mathbb{S}^3$.
In this way, we recover the tessellation of $\mathbb{S}^3$ by seven dumplings.
\par\smallskip
To make our idea and a possible inductive process clear, we start this chapter by studying a chamber bounded by five (as well as six)  
hyperplanes stable under the action of the cyclic group $\mathbb{Z}_5$ (resp. $\mathbb{Z}_6$); 
though such a chamber is not unique, the above statement with obvious modification is still true, 
if we understand the special edge of $D_i$ as $D_{i-1}\cap D_i\cap D_{i+1}$. A 3-dumpling is a tetrahedron, and a 4-dumpling is a prism. After the collapsing process as above, we get the manifolds $X_5$ (resp. $X_6$). We prove in detail that $X_5$ is homeomorphic to $\mathbb{S}^3$; for other cases proof is similar.  
\par\smallskip
\subsection{Five hyperplanes in the 4-space}
Five hyperplanes in the projective 4-space cut out sixteen 4-simplices. 
If a hyperplane is at infinity, then the remaining four can be considered as the four coordinate hyperplanes, 
dividing the space into $2^4$ chambers. The boundary of each 4-simplex is of course made by five tetrahedra. 
In other words, five tetrahedra are glued together to make a (topological) 3-sphere: around a tetrahedron, 
glue four tetrahedra along their faces, and further use the valley of each pair of triangles in order to close the pair like a book. 
\subsubsection{Labeling and a study of 2-skeleton}\label{45}
We study the above tessellation more precisely. Label the five hyperplanes as $H_1,\dots,H_5$.
Choose one 4-simplex, out of sixteen, and call it CC. We use the convention
$$H_{ij}=H_i\cap H_j,\quad  H_{ijk}=H_i\cap H_j\cap H_k, \ \dots$$
The five vertices of CC are:
$$H_{2345},\quad H_{1345},\ \cdots,\ H_{1234}$$
which are often denoted simply by 
$$ 1,\ 2,\ \dots,\ 5,$$
respectively. The boundary of CC is made by five tetrahedra: 
$$D_1=2345,\quad D_2=1345,\ \dots.$$
(Tetrahedron $D_1$ has vertices $2,3,4,5$.) Any two of them share a triangle:
$$D_i\cap D_j=\Delta(\{1,\dots,5\}-\{i,j\}),\quad i\ne j;$$
there are ten of them.
These tetrahedra tessellate the 3-sphere $\mathbb{S}^3$. We consider the sequence $D_1,D_2,\dots$, in this order, and name the intersections of the two consecutive ones as:
$$D_1\cap D_2=\Delta(345),\quad D_2\cap D_3=\Delta(451),\ \cdots D_5\cap D_1=\Delta(234).$$
(Note that the union of these five triangles form a M\"obius strip, but for a while just forget it.)
Remove these five triangles from the 2-skeleton of the tessellation. Then five triangles remain, forming another M\"obius strip
$$M:\quad
\scriptsize\begin{array}{ccccccc  }
3&-&4&-&5&-&1\\
 \quad\ba& & /\quad\ba &&/\quad\ba && /\quad\\
    &1&-&2&-&3& 
\end{array}.$$\normalsize
Its boundary 
$$C=\partial M:\quad 1 \ -\  2 \ -\  3 \ -\  4 \ -\  5 \ -\  1$$
is an unknotted circle because $C$ bounds also the disc
$$\Delta(234)\cup  \Delta(124)\cup \Delta(451): \qquad
\scriptsize\begin{array}{ccccc  }
3&-&4&-&5\\
 \quad\ba& & /\quad\ba &&/\quad\\
    &2&-&1&
\end{array}.$$\normalsize

\subsubsection{Visualization}
\begin{figure}[h]
\begin{center} 
\includegraphics[width=140mm]{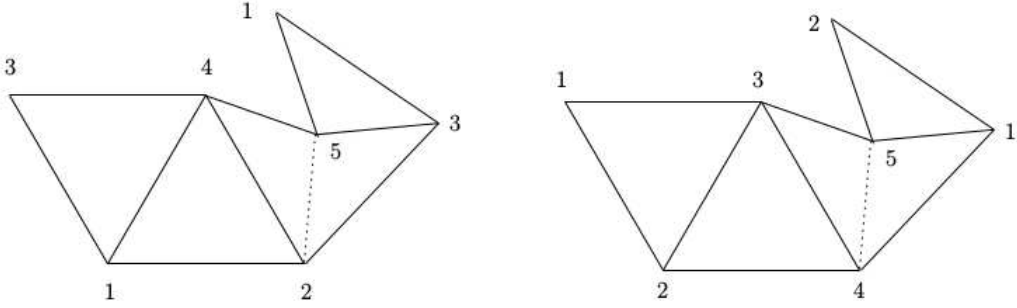}
\end{center}
\caption{M\"obius strips $M$ (left) and the complementary one (right) before identification along appropriate edges.}
\label{fig8}
\end{figure}

Since the situation is quite simple, we visualize the happenings. Consider a  tetrahedron with vertices $\{1,2,3,4\}$ in our space, and  
put vertex 5 at the barycenter. 
The barycentric subdivision of this tetrahedron (which yields $D_1,\dots,D_4$) together with its complementary tetrahedron $D_5$ fills the 3-sphere. 
Figure \ref{fig8} shows the 2-skeleton consisting of ten triangles as the union of the M\"obius strips $M$ and the complementary one. 
You can find the disc, given in the previous subsection,  bounded by the curve $C$. When you make a paper model of the 2-skelton out of Figure \ref{fig8}, be ware that the two M\"obius strips have different orientations.

\begin{figure}
\begin{center}
\includegraphics[width=100mm]{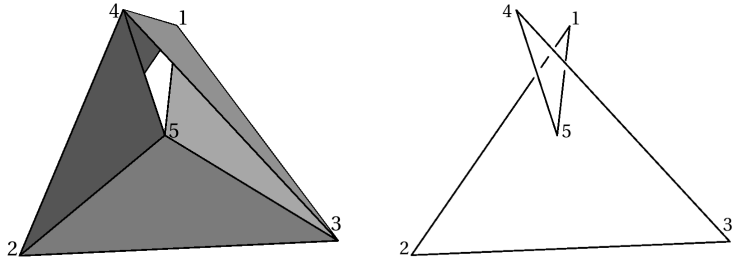}
\end{center}
\caption{(left) M\"obius strip $M$. (right) Boundary $\partial M$.}
\label{fig7}
\end{figure}

\subsubsection{From a solid torus}
We would like to give another description, which will prepare the next step where the number $m$ of hyperplanes is $m=6,7,\dots$, 
in such cases, naive description above would hardly work. We start with a {\it vertical triangular prism}, 
horizontally sliced into five  thinner triangular prisms. Its top and the bottom triangles are identified after a $4\pi/3$ rotation; 
this makes the vertical prism a solid torus. The boundary torus is tessellated by fifteen rectangles. Now collapse vertically five diagonal rectangles. 
This changes every remaining rectangles into triangles, and the thinner triangular prisms into tetrahedra. 
Further, use each collapsed edge as a valley of each pair of triangles face-to-face in order to close the pair like a book. 
This identification makes the solid torus a 3-sphere tessellated by five tetrahedra, recovering the tessellation described in \S\ref{45}.
\par\medskip\noindent{\bf A vertical triangular prism horizontally sliced.}
We show/repeat the above process by using the labels $1,2,3,4,5$ (they correspond the labels of the five hyperplanes $H_1,\ H_2,\dots$). 
Consider five triangles labeled by two consecutive numbers: $51,\ 12, \dots$ (they correspond the triangle $H_{51}\cap{\rm CC},\dots$). 
The triangle 51, for example, has three vertices with labels 
$$51:\quad\{5124,\ 5123, \ 5134\}$$ (they correspond the points $\{H_{4512},H_{5123},H_{3451}\}$). 
Now operate the group $\mathbb{Z}_5$ on the labels: $j\to j+1$ mod 5, and we get four other triangles, for example, the triangle 12 has three vertices:
$$12:\quad\{1235, \ 1234, \ 1245\}.$$
We make a triangular prism by putting the triangle 51 above the triangle 12;
the two vertices of 51 having the same labels of the ones of 12 should be put just above each other. 
In a word this prism is the sandwich made by the two triangles 51 and 12.

\begin{figure}
\begin{center}
\includegraphics[width=140mm]{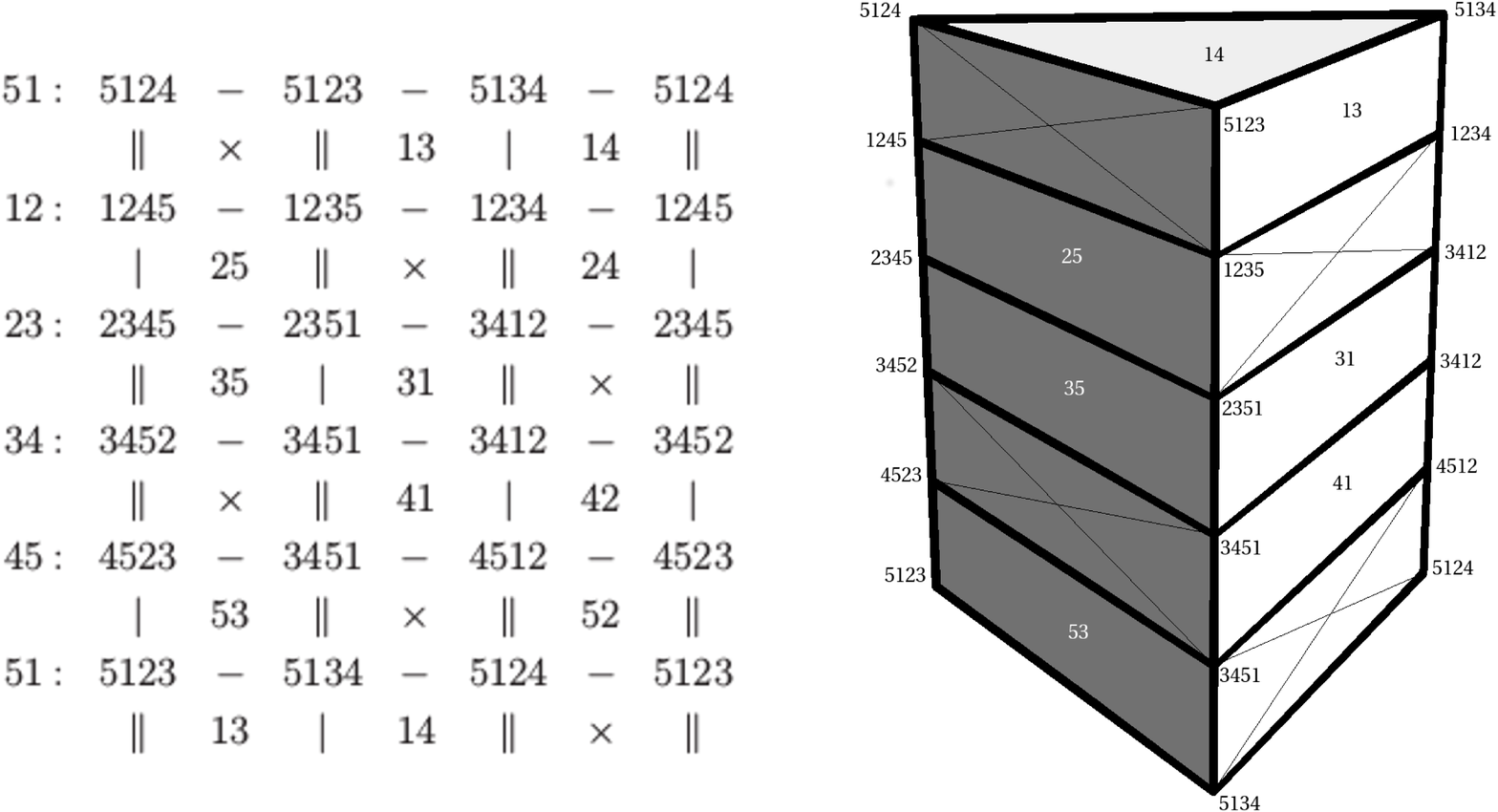}
\end{center}
\caption{(left) Table of five piles. (right) Triangular prism split up into five slices. 
The label 14 denotes the upper back rectangle.}
\label{fig1}
\end{figure}

\noindent In this way, we make five sandwiches and pile them to form a vertical triangular prism. 
The table in Figure \ref{fig1} repeats what we described. Each off-diagonal rectangle is labeled by the two numbers common to the four vertices.

\par
\begin{figure}
\begin{center}
\includegraphics[width=120mm]{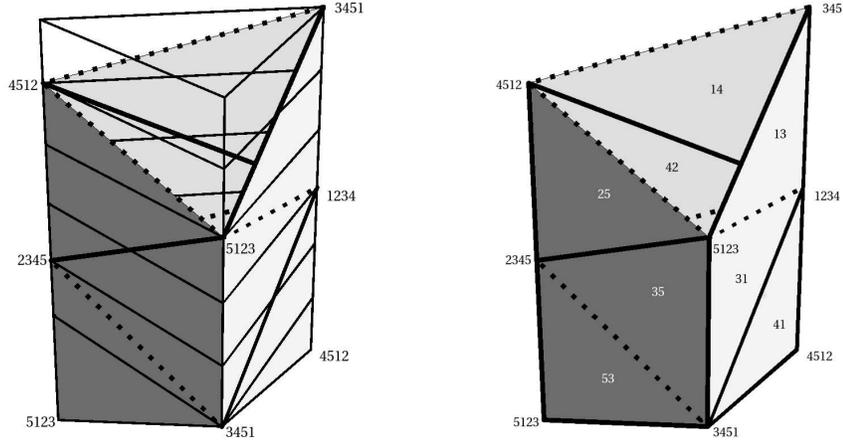}
\end{center}
\caption{(left) Collapsing rectangles. (right) Tiling the boundary of the solid cylinder by ten triangles. 
The dashed line is the special loop $C$.}
\label{fig2}
\end{figure}

\medskip\noindent{\bf Collapsing rectangles.}
If we identify two points with the same indices {\it as sets}, for example, $4512=5124=1245$, 
you find five rectangles with pairwise coincide vertices (marked by $\times$), situated diagonally. 
We collapse such rectangles vertically to a segment, each called a {\bf special edge}. 
Note that a special edge has two ends which are labeled by consecutive numbers: $\{1234,2345\},\{2345,3451\},\dots$ 
Accordingly all the remaining rectangles become triangles and, the five thin prisms become tetrahedra. 
But as a whole it still remains to be a cylinder. This collapsing process is described in Figure \ref{fig2}. Here the prisms are hollow so that we can peep inside.

\par\medskip\noindent{\bf Identifying the top and the bottom to get a solid torus $ST$.} 
We identify the top triangle and the bottom triangle of the collapsed vertical prism according to the labels; note that we must twist the triangle by $4\pi/3$. 
By the identification we get a solid torus, say $ST$; its meridean ${\bf mer}$ is represented by the boundary of a(ny) triangle, in Figure \ref{fig1}, 
say 51 (or 12, 23,\dots). Note that it is still an abstract object. The five special edges now form a circle:
$$C:\quad 1234-2345-3451-4512-5123-1234.$$
We take a vertical line, in Figure \ref{fig1} (left), with a $4\pi/3$ twist, say,
$${\bf par}: 5123 - 3451 - 5123=4-2-4$$
as a parallel.
\par
\begin{figure}
\begin{center}
\includegraphics[width=120mm]{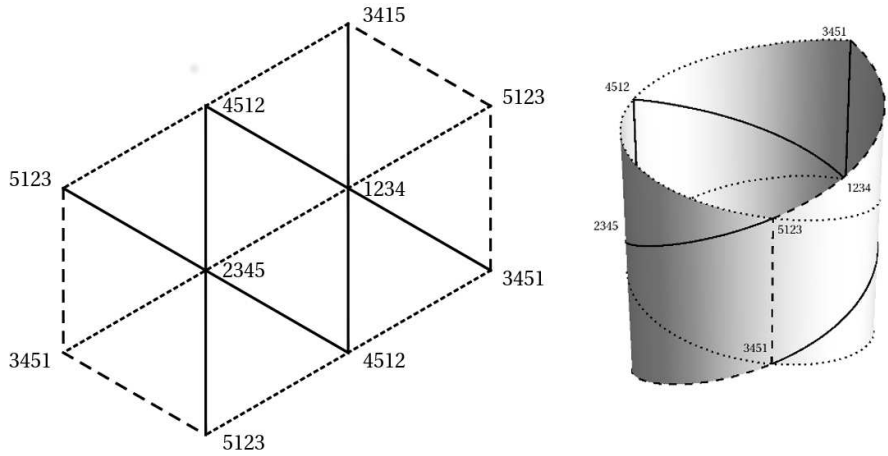}
\end{center}
\caption{(left) Ten triangles tesselating the hexagon. (right) Identifying the left and right vertical sides we get a cylinder topologically
equivalent to that of Figure \ref{fig2}-right.
The special curve $C$ is dotted while $ST$-parallel is dashed.}
\label{fig3}
\end{figure}

\medskip\noindent{\bf The torus $T$ and the solid torus $UT$.}
We consider in our space ($\sim \mathbb{S}^3$) a usual (unknotted) solid torus, say $UT$, and put our solid torus $ST$ fills outside of $UT$, 
so that the tessellated boundary torus
$$T:=\partial ST\ (=\partial UT)$$
has {\bf par} as the meridean, and {\bf mer} as the parallel.
From the   description of  the vertical prism, we see ten triangles tessellating the torus $T$ arranged in a hexagonal way (Figure \ref{fig3}-left). 

\begin{figure}
\begin{center}
\includegraphics[width=100mm]{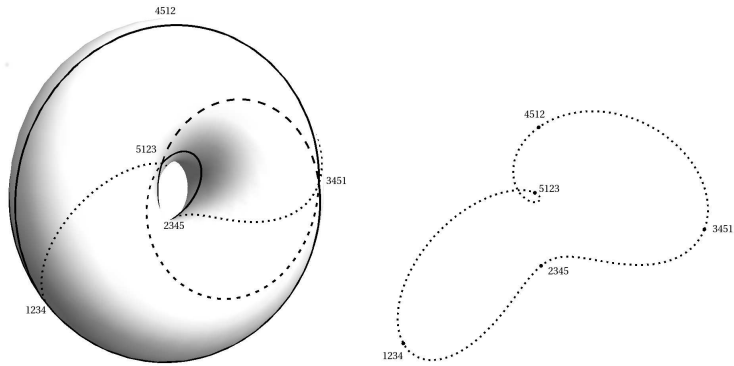}
\end{center}
\caption{(left) Tesselation of the torus $\partial ST$. Special curve $C$ is dotted while $ST$-parallel is dashed. (right) The special curve $C$.}
\label{fig4}
\end{figure}

Note that the triangulated torus above is not a simplicial complex; indeed there are two edges with the same vertices, for example, $5123-3451$.
Anyway, you can see the curve $C$ lies on the torus $T=\partial ST$ as a $(2,1)$-curve (see Figure \ref{fig4}).

\par\medskip\noindent{\bf Collapsing the solid torus $UT$ by folding $T$.} 
On the torus $T$,  there are pairs of triangles with the same vertices. We identify these triangles  by folding the torus $T$ along the curve $C$.
This collapses the solid torus $UT$; we prove below that the resulting space, say $X_5$, is still $\mathbb{S}^3$.
\begin{figure}
\begin{center}
\includegraphics[width=100mm]{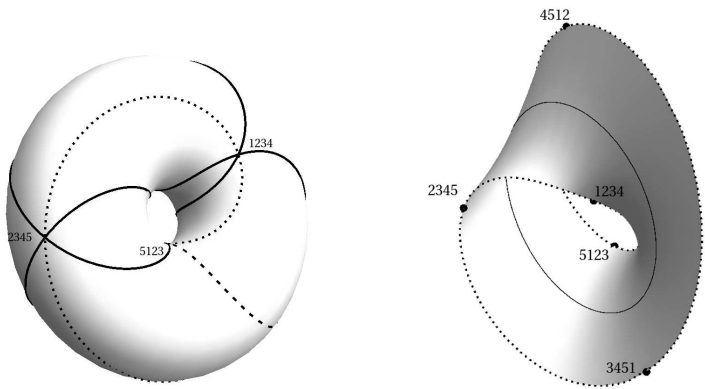}
\end{center}
\caption{(left) Tesselation of the torus $\partial UT$. Special curve $C$ is dotted while $UT$-meridian is dashed. (right) M\"obius strip $M$.}
\label{fig5}
\end{figure} 
The torus $T$ is folded to be a M\"obius strip $M$ with boundary $C$:
\scriptsize
$$M:\quad\begin{array}{ccccccc  }
4512&-&5123&-&1234&-&2345\\
 \quad\ba&25 & /\ 35\ \ba &31&/\ 41\ \ba &42& /\quad\\
    &2345&-&3451&-&4512& 
\end{array}$$
\normalsize
\begin{figure}
\begin{center}
\includegraphics[width=120mm]{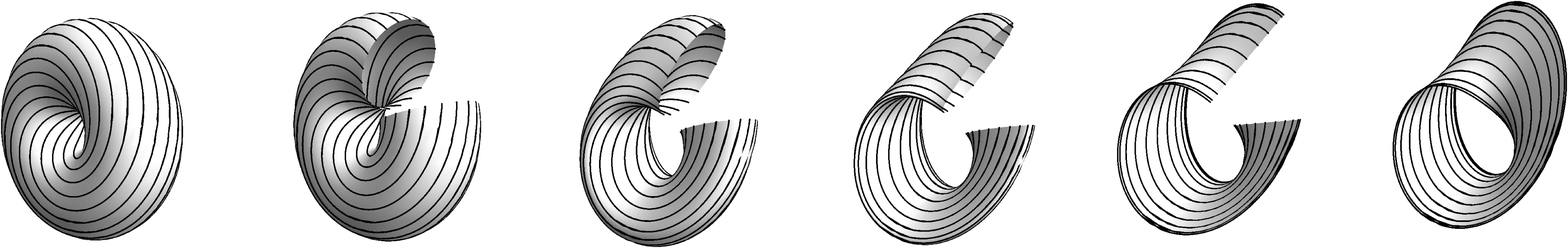}
\end{center}
\caption{(from left to right) Folding the torus $T$ along the special curve $C$ to be the M\"obius strip $M$. A sector is cut out to show the folding.}
\label{fig6}
\end{figure} 
Note that $M$ is a simplicial complex. If we use complementary labeling, for example, $5$ for $1234$, this is exactly the same as we got several pages before. 
On the hexagonal expression of the torus $T$, this folding is done as closing a book along the diagonal.

\begin{figure}
\begin{center}
\includegraphics[width=60mm]{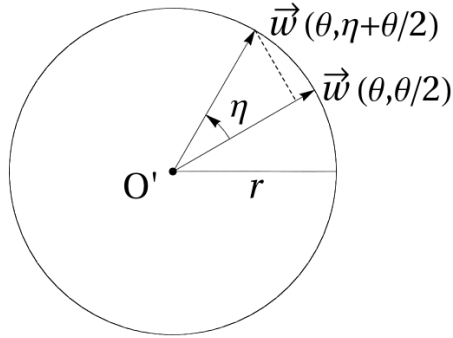}
\end{center}
\caption{Meridian section of $UT$. The homotopy collapsing $UT$ onto $M$ follows the dashed line.}
\label{fig10}
\end{figure}

The process of collapsing the solid torus $UT$ onto 
the M\"obius strip $M$ can be explicited by the following homotopy parametrized by $t\in[0,1]$:
$$
F(\vartheta,\eta,r,t)=O'+(1-t)r\overrightarrow w(\vartheta,\eta+\frac{\vartheta}{2})+tr\cos\eta\cdot\overrightarrow w(\vartheta,\frac{\vartheta}{2}),
$$
where 
$$
O'=\sqrt{2}(\cos\vartheta,\sin\vartheta,0),\quad \overrightarrow w(\vartheta,\eta)=(\cos\eta\cos\vartheta,\cos\eta\sin\vartheta,\sin\eta)
$$
and $(\exp(i\vartheta),\exp(i\eta),r)$ parametrizes the solid torus $\mathbb{S}^1\times\mathbb{D}^2$ (Figure \ref{fig10}).

\par\medskip\noindent{\bf Puffing the M\"obius strip to make a torus.} The inverse operation of folding the (unknotted) torus $T$ 
to the M\"obius strip $M$ (with $\pm\pi$ twist) can be described as follows: Consider a M\"obius strip $M$ being double sheeted like a flat tire folded along $C$. 
Then puff the tire to be a air-filled tire $T$. It is illustrated in Figure \ref{fig6} seen from right to left. 
\par\medskip\noindent{\bf A disc bounding $C$.} Let us repeat the disc in \S \ref{45} bounding $C$, using complementary labeling:
\scriptsize
$$\begin{array}{ccccc  }
4512&-&5123&-&1234\\[2mm]
 \quad\ba&51 & /\quad35\quad\ba &23&/\quad\\[2mm]
    &3451&-&2345&
\end{array}.$$\normalsize

\par\medskip\noindent{\bf Proof that $X_5$ is homeomorphic to $\mathbb{S}^3$.} 
\noindent The collapsing of $\mathbb{S}^3=UT\cup_T ST$ by folding $T$ along $C$ and collapsing $UT$ on $M$ is nothing but the gluing 
of the set of the five tetrahedra tiling $ST$ in such a way that the tetrahedra facets are identified by pairs. In doing so we must get $ST$ back in one hand,
and, on the other hand, the triangles of the boundary $\partial ST$ are identified by pairs to get the M\"obius strip $M$. 
The result $X_5$, is a closed $3$-dimensional manifold since the link of each vertex is homeomorphic to a $2$-sphere, namely the $2$-skeleton of a tetrahedron. 
For instance the link of the vertex $1234$ is the union of the four triangles $\{15,25,35,45\}$, where 15 for example, is the triangle with vertices $1523,1534$ and $1542$. 
The $3$-manifold $X_5$, as a continuous image 
of $\mathbb{S}^3$, is clearly connected.

\begin{figure}
\begin{center}
\includegraphics[width=90mm]{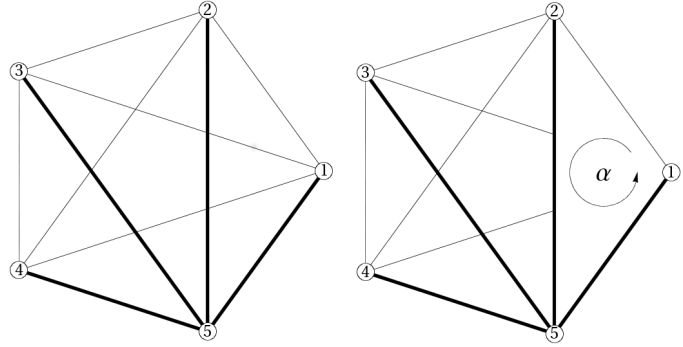}
\end{center}
\caption{(left) The complete graph $X_5^1$ over five vertices and a maximal tree base at $5$. (right) One of the six generators of $\pi(X_5^1,5)$:
$\alpha=(5,1,2)$.}
\label{fig9}
\end{figure}

\noindent Now, computing the fundamental group $\pi(X_5,\ast)$ will prove that $X_5$ is simply connected and therefore, by Poincar\'e-Perelman theorem, 
homeomorphic to $\mathbb{S}^3$. Indeed, since $X_5$ is endowed with a structure of simplicial complex, $\pi(X_5,\ast)$ is isomorphic to $\pi(X_5^2,\ast)$, 
where $X_5^2$ denotes the $2$-skeleton. The group $\pi(X_5^2,\ast)$ can be described by generators and relators as follows. 
Firstly, we observe that the $1$-skeleton $X_5^1$ is the complete graph over the five vertices $1,2,3,4,5$ (using complementary labeling) (Figure \ref{fig9}). 

Recall that we can use a maximal tree (that is a tree containing all vertices) to compute the fundamental group $\pi(X_5^1,5)$. 
Indeed, to each egde not contained in the tree, we associate in a obvious way a loop based at $5$. 
The set of these loops is a basis of the free group $\pi(X_5^1,5)$.
Therefore we see that
$\pi(X_5^1,\ast)$ is the free group over the six following generators (still using complementary labeling for the sequence of vertices 
and chosing $5$ as a base point):
$$
\alpha=(5,1,2),\quad\beta=(5,1,3),\quad\gamma=(5,1,4),\quad\delta=(5,2,3),\quad\varepsilon=(5,2,4),\quad\zeta=(5,3,4)
$$
Secondly, the relators of $\pi(X_5^2,\ast)$ are associated with the ten facets coming from the five tetrahedra tiling $X_5$: 
$$
12,\quad13,\quad14,\quad15,\quad23,\quad24,\quad25,\quad34,\quad35,\quad45.
$$
\begin{figure}
\begin{center}
\includegraphics[width=100mm]{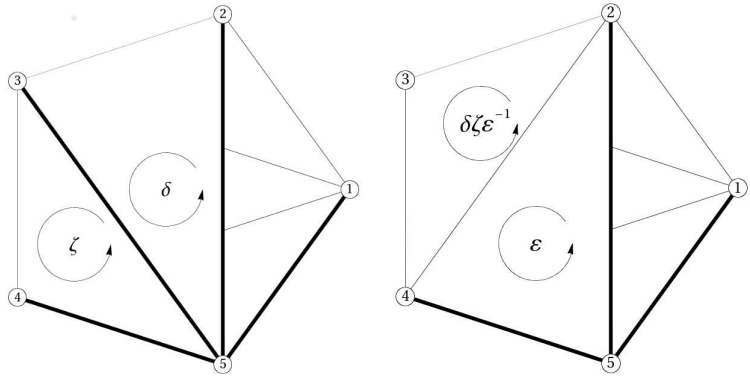}
\end{center}
\caption{(left) The generators $\delta=(5,2,3)$ and $\zeta=(5,3,4)$. (right) The generator $\varepsilon=(5,2,4)$ 
and the relator $\delta\zeta\varepsilon^{-1}=(2,3,4)$.}
\label{fig11}
\end{figure} 
It yields to the following relators (Figure \ref{fig11}):
$$
\zeta,\quad\varepsilon,\quad\delta,\quad\delta\zeta\varepsilon^{-1},\quad\gamma,\quad\beta,\quad\beta\zeta\gamma^{-1},\quad\alpha,
\quad\alpha\varepsilon\gamma^{-1},\quad\alpha\delta\beta^{-1},
$$
something which proves that 
$$
\begin{array}{l}
 \pi(X_5,\ast)=\pi(X_5^2,\ast)=\\
<\alpha,\beta,\gamma,\delta,\varepsilon,\zeta:
\zeta=\varepsilon=\delta=\delta\zeta\varepsilon^{-1}=\gamma=\beta=\beta\zeta\gamma^{-1}=\alpha=\alpha\varepsilon\gamma^{-1}=\alpha\delta\beta^{-1}=1>\\
=1.
\end{array}
$$

\par\medskip\noindent{\bf Conclusion.} 
\noindent
In the sphere $\mathbb{S}^3=UT\cup_T ST$, the solid torus $UT$ is collapsed to an unknotted M\"obius strip $M$, 
which is twisted only by $\pm\pi$ (see Figure \ref{fig6}), and the resulting $3$-manifold $X_5$ is a $3$-sphere.
In this way we recover the tessellation of $\partial \Delta_4=\partial {\rm CC}$ described in \S \ref{45}.


\subsection{Six hyperplanes in the 4-space}
Six hyperplanes in the projective 4-space cut out six 4-simplices, fifteen prisms of type $\Delta_3\times\Delta_1$ and ten prisms of type $\Delta_2\times\Delta_2$. If a hyperplane is at infinity, then the remaining five bound a simplex. Other chambers touch this simplex along 3-simplices (five $\Delta_3\times\Delta_1$ ), along 2-simplices (ten $\Delta_2\times\Delta_2$), along 1-simplices (ten $\Delta_3\times\Delta_1$) and along 0-simplex (five $\Delta_4$). 

Though a prism of type $\Delta_3\times\Delta_1$ is bounded by six hyperplanes, it does not admit an action of the group $\mathbb{Z}_6$. So we consider one of the prisms of type $\Delta_2\times\Delta_2$, which is bounded by six prisms (of type $\Delta_2\times\Delta_1$). We study how these six prisms tessellate a 3-sphere.
\subsubsection{Labeling and a study of 2-skeleton}
Let the six hyperplanes $H_1,\dots,H_6$ bound a chamber CC of type $\Delta_2\times\Delta_2$. The chamber CC has six faces $D_1,\dots,D_6$, which form the $\mathbb{Z}_6$-orbit of the prism
$$D_1={\rm CC}\cap H_1:\scriptsize\begin{array}{ccccc}
       1346& --&--&--            &1456\\
         |& \ba\quad&4&\quad/                &|\\
         |& 3\quad1234&--&  1245\quad5      &|\\
         |& /\quad&   2     &\quad\ba       &|\\
       1236& --&--&--           &1256\end{array}\normalsize,$$
where $1234=H_{1234}=H_1\cap\cdots\cap H_4$. Note that the rectangular face of $D_1$ behind is $D_1\cap D_6$, and the rectangular face in front below is $D_1\cap D_2$.
\footnote{Using the notation introduced in Chapter 2, CC is represented by $++++++$, and its boundary consistes of 
$$D_1=0+++++,\quad D_2=-0++++,\ \dots,\ D6=-----0.$$
The boundary of $D_1$ consists of (cf. \S \ref{2D})
$$\begin{array}{lll}
00++++:{\rm4-gon}&
0-0+++:{\rm triangle}&
0--0++:{\rm rectangle}\\
0---0+:{\rm triangle}&
0----0:{\rm4-gon}
\end{array}\normalsize
$$} 
The 2-skeleton of CC consists of six triangles and nine rectangles. Remove the intersections  of the two consecutive ones:
$$D_1\cap D_2,\quad D_2\cap D_3,\quad \dots,\quad D_6\cap D_1:$$
$$\scriptsize\begin{array}{cccccccccc}
5612& -&6123& -&2356&  &    &  &&\\
  | &12&  | &23& |  &  &    &  &&\\
1245& -&1234& -&2345& -&4512&  &&\\
    &  &  | &34& |  &45& |  &  &&\\
    &  &3461& -&3456& -&4561& -&6134&\\
    &  &    &  & |  &56& |  &61& |  &\\
    &  &    &  &5623& -&5612& -&6123&.
  \end{array}\normalsize$$
  (Note that the union of these six rectangles form a cylinder, {\it not} a M\"obius strip.)
Then  the remaining six triangles and three rectangles form a M\"obius strip
$$M:\quad\scriptsize\begin{array}{ccccccccc  }
5612&-&6123&-&1234&-&2345&-&3456\\
 \quad\ba&26 & /\quad\ba &31&/\quad\ba  &42&/\quad\ba&53&/\\
    &2356& &3461& &1245& &2356\\
 &\quad\ba&36 & /\ 46\ \ba &41&/\ 51\ \ba  &52&/\quad   & \\
 &   &5634&-&6145&- &5612&&\\
\end{array}\normalsize$$
\begin{center} M\"obius  6
\end{center}
Its boundary 
$$C=\partial M:\quad 1234-2345-3456-4561-5612-6123-1234$$
is an unknotted circle because $C$ bounds also the disc
$$\scriptsize\begin{array}{ccccccc}
&6123\quad&-&\quad1234&\\[2mm]
\quad/&&\ba\quad13\quad/&&\ba\quad\\[2mm]
5612&61\quad&6134&34&2345\\[2mm]
\quad\ba&&/\quad46\quad\ba&&/\quad\\[2mm]
&4561\quad&-&\quad3456&
\end{array}\normalsize.$$

\subsubsection{From a solid torus}
We start with a vertical rectangular prism, horizontally sliced into six thinner rectangular prisms. Its top and the bottom rectangles are identified after a $4\pi/4$ rotation; this makes the vertical prism a solid torus. The boundary torus is tessellated by 24 rectangles. Now collapse vertically six diagonal rectangles. This changes eighteen rectangles (out of 24) into six segments (forming a circle) and twelve triangles. Accordingly, the six thinner rectangular prisms (actually cubes) into prisms. Next use each collapsed edge as a valley of each pair of triangles face-to-face in order to close the pair like a book, and further identify three pairs of rectangles. This identification changes the solid torus into a 3-sphere tessellated by six prisms.

\par\medskip\noindent{\bf A vertical rectangular prism horizontally sliced.}
We show/repeat the above process by using the labels $1,2,3,4,5,6$ (they correspond the labels of the six hyperplanes $H_1,\ H_2,\dots$). Consider six rectangles labeled by two consecutive numbers: $61,\ 12, \dots$ (they correspond the rectangle $H_{61}\cap{\rm CC},\dots$). The rectangle 61, for example, has four vertices with labels 
$$61:\quad\{6125,\ 6123, \ 6134,\ 6145\}$$ (they correspond the points $\{H_{5612},H_{6123},H_{6134},H_{6145}\}$). Now operate the group $\mathbb{Z}_6$ on the labels: $j\to j+1$ mod 6, and we get five other rectangles, for example, the rectangle 12 has three vertices:
$$12:\quad\{1236, \ 1234, \ 1245,\ 1256\}.$$
We make a rectangular prism by putting the rectangle 61 above the rectangle 12;
the two vertices of 61 having the same labels of the ones of 12 should be put just above each other. In a word this prism is the sandwich made by the two rectangles 61 and 12.

In this way, we make six sandwiches and pile them to form a vertical rectangular prism. The following table repeats what we described. Each off-diagonal rectangle is labeled by the two numbers common to the four vertices.

$$\scriptsize\begin{array}{rccccccccc}
61:&6125&-&6123&-&6134&-&6145&-&6125\\
  &\|   &\times&\|  &13&|   &14&|   &15&\|   \\
12:&1256&-&1236&-&1234&-&1245&-&1256\\
  &|   &26  &\| &\times&\|   &24&|   &25&|  \\
23:&2356&-&2361&-&2341&-&2345&-&2356\\
  &|   &36  &|  &31&\| &\times&\|   &35&|   \\
34:&3456&-&3461&-&3412&-&3452&-&3456\\
  &\|   &46  &|  &41&| &42&\|   &\times&\|   \\
45:&4563&-&4561&-&4512&-&4523&-&4563\\
  &\|   &\times&\|  &51&| &52&|   &53&\|   \\
56:&5634&-&5614&-&5612&-&5623&-&5634\\
  &|   &64&\|  &\times&\| &62  &|   &63&|   \\
61:&6134&-&6145&-&6125&-&6123&-&6134\\
  &|   &14&|  &15&\| & \times &\|   &13&|   \\
\end{array}\normalsize$$
\begin{center} Table of piles 6
\end{center}
\par\medskip\noindent{\bf Collapsing rectangles.}
If we identify two points with the same {\it consecutive} indices {\it as sets}, for exmaple, $6123=1236=2361$ but $2356\not=5623$, you find six rectangles with pairwise coincide vertices (marked by $\times$), situated diagonally. We collapse such rectangles vertically to a segment, each called a {\bf special edge}. Note that a special edge has two ends which are labeled by consecutive numbers: $\{1234,2345\},\{2345,3451\},\dots$ Accordingly the remaining rectangles next to the special edges become triangles, and the six thin rectangular prisms become triangular prisms ($\sim \Delta_2\times\Delta_1$). But as a whole it still remains to be a cylinder. 
\par\medskip\noindent{\bf Identifying the top and the bottom to get a solid torus $ST$.} 
We identify the top rectangle and the bottom rectangle of the collapsed vertical prism according to the labels; note that we must twist the triangle by $4\pi/4$. By the identification we get a solid torus, say $ST$; its meridean ${\bf mer}$ is represented by the boundary of a(ny) rectangle, in Table of piles 6, say 61 (or 12, 23,\dots). Note that it is still an abstract object. The six special edges now form a circle:
$$C: 1234-2345-3456-4561-5612-6123-1234.$$
We take a vertical line, in Table of piles 6, with a $4\pi/4$ twist, say,
$${\bf par}: 6123 - 3461 - 4561 - 6123$$
as a parallel.

\par\medskip\noindent{\bf The torus $T$ and the solid torus $UT$.}
We consider in our space ($\sim \mathbb{S}^3$) a usual (unknotted) solid torus, say $UT$, and think our solid torus $ST$ fills outside of $UT$, so that the tessellated boundary torus
$$T:=\partial\ ST\ (=\partial\ UT)$$
has {\bf par} as the meridean, and {\bf mer} as the parallel.
From the   description of  the vertical prism, we see twelve triangles and six rectangles tessellating the torus $T$ arranged in a hexagonal way (identify the opposite sides of the hexagon):
$$\scriptsize\begin{array}{ccccccccccc}
5612&-&5623&-&3456&&& & &\\
|   &26\ba62&|  &63&|   &64\ba&  &&   &\\
2356&-&6123&-&6134&-&4561&&&\\
|   &36 &| &31\ba13&|   &14 &|   &15\ba &  &\\
3456&-&3461&-&1234&-&1245&-&5612&\\
  & \ba46 &|  &41 &| &42\ba24&|   &25 &|   &\\
    & &4561&-&4512&-&2345&-&2356&\\
    & &  &\ba51 &| &52 &|   &53\ba35&|   & \\
    &&&&5612&-&5623&-&3456&&\\
  \end{array}\normalsize$$
\begin{center} Torus (hexagon) 6
\end{center}
\par\medskip\noindent{\bf Collapsing the solid torus $UT$ by folding $T$.}
We identify two vertices with the same indices as sets. Along each collapsed edge (special edge), there are two triangles with the same vertices; these triangles are folded along the special edge into one triangle. Further there are three pairs of rectangles with the same vertices; they are also identified. Consequently the torus $T$ is folded along the curve $C$ to be a M\"obius strip, which turns out to be the same as the one in the previous subsection (shown as M\"obius 6).
\comment{
$$M:\quad\scriptsize\begin{array}{ccccccccc  }
5612&-&6123&-&1234&-&2345&-&3456\\
 \quad\ba&26 & /\quad\ba &31&/\quad\ba  &42&/\quad\ba&53&/\\
    &2356& &3461& &1245& &2356\\
 &\quad\ba&36 & /\ 46\ \ba &41&/\ 51\ \ba  &52&/\quad   & \\
 &   &5634&-&6145&- &5612&&\\
\end{array}\normalsize$$
\begin{center} M\"obius  6
\end{center}}
In this way, we get a space $X_6$ from the solid torus $ST$ made by six prisms by folding the torus $T$ (collapsing $UT$).

Since the curve $C$ is unknotted, this folding is done exactly the same as in the previous section. Then $X_6$ is homeomorphic to $\mathbb{S}^3$.

\par\medskip\noindent{\bf A model.} 
Put three vertices with labels without numeral 6: $1245, 2345, 1234$ inside the prism $(6123,5612,6235)*(6134,6145,6345)$. 
The M\"obius strip above is shown in Figure \ref{Moe6} (\cite{CYY}). From this picture, one can hardly see the $\mathbb{Z}_6$ action. 
\begin{figure}[h]
\begin{center} 
\include{Moe6s} 
\end{center}
\caption{Five prisms pack a prism. M\"obius strip made by six triangles and three rectangles is shown(Moe6)}
\label{Moe6}
\end{figure} 

\subsubsection{$6-1=5$}
If we remove the sixth hyperplane, then the remaining five hyperplanes bounds a 4-simplex. This process can be described as follows: the prism on the sixth plane reduces to a segment. More precisely, the two triangles of the boundary of the sixth prism reduce to two points, and the prism reduces to a segment connecting these two points. Combinatorial explanation: The prism on the sixth plane with vertices
$$\scriptsize\begin{array}{ccccc}
1256&-&-&-&1456\\
|&\ba\qquad&&\qquad/&|\\
|&2356&-&3456&|\\
|&/\qquad&&\qquad\ba&|\\
1236&-&-&-&1346\end{array}\normalsize\quad {\rm on\ 6}$$
reduces ($1256,2356,1236\to1235,\ 1456,3456,1346\to1345$) to the segment
$${\bf 1235\ ---\ 1345}.$$
Note that 
$$\begin{array}{ll}
\{1,2,3,5\}&=\{1,2,5,6\}\cup\{2,3,5,6\}\cup\{1,2,3,6\}-\{6\},\\[2mm]
\{1,3,4,5\}&=\{1,4,5,6\}\cup\{3,4,5,6\}\cup\{1,3,4,6\}-\{6\}.
\end{array}$$
Accordingly the other prisms reduce to tetrahedra, for example:
$$ \begin{array}{cccl}
\scriptsize\begin{array}{ccc}
3412&---&3612\\
|&\ba\qquad\qquad/&|\\
|&4512-5612&|\\
|&/\qquad\qquad\ba&|\\
3452&---&3562\end{array}\normalsize
&\longrightarrow 
&\scriptsize\begin{array}{cl}
3412&\\
|&\ba\qquad\ \ba\\
|&4512-{\bf1235}\\
|&/\qquad\ /\\
3452&\end{array}\normalsize
 &\normalsize{\rm on\ 2} \ ({\rm similar\ on}\ 4),\\
\\
\scriptsize\begin{array}{ccc}
6145&---&6345\\
|&\ba\qquad\qquad/&|\\
|&1245-2345&|\\
|&/\qquad\qquad\ba&|\\
6125&---&6235\end{array}\normalsize&\longrightarrow 
&\scriptsize\begin{array}{c}
{\bf1345}\\
/\quad |\quad \ba\\
1245-2345\\
\ba\quad |\quad /\\
{\bf1235}\end{array}\normalsize &{\rm on\ 5}\ ({\rm similar\ on}\ 1\ {\rm and}\ 3).
 \end{array}\normalsize$$

\par\noindent
{\bf Remark}. Consider an $n$-gon in the plane bounded by $n$ lines. If you remove the $n$-th line, then you get an $(n-1)$-gon. During this process, the $n-3$ sides away from the $n$-th sides do not change. Consider an $n$-hedron in the 3-space bounded by $n$ planes. If you remove the $n$-th plane, you get an $(n-1)$-hedron. During this process, faces away from the $n$-th face remain unchanged. In the present case, all the prism-faces change into tetrahedra.

\subsubsection{$5+1=6$\label{56}}
If we add the sixth hyperplane to the five hyperplanes, then among the sixteen 4-simplices $\D_4$, some are untouched, some are divided into $\D_4$ and $\D_3\times\D_1$ (with section $\D_3$), and some are divided into $\D_2\times\D_2$ and $\D_3\times\D_1$ (with section $\D_2\times \D_1$.)
We will describe the change of $\D_4$ into $\D_2\times\D_2$ via that of these boundaries: the 3-sphere tessellated by five tetrahedra into that tessellated by six prisms. 
\par\smallskip
As a preparation, we would like to explain this cutting process by using a 1-dimension less model: we truncate a vertex of a polyhedron in 3-space. 
Let us consider a polyhedron in the 3-space, say for example a tetrahedron bounded by four faces $\{0,1,2,3\}$, 
and cut the vertex $123=1\cap2\cap3$ by the new plane named 4. The tetrahedron is truncated and is divided into a prism and a small tetrahedron; 
they share a triangle which should be called 4. Instead of describing this process by 3-dimensional pictures, 
we express it by 2-dimensional pictures describing the change near the vertex $123$: the plane is divided into three chambers $1,2,3$ 
(imagine the letter Y) sharing a unique point representing the vertex 123, which we call the center. 
The cutting is expressed by a blowing-up at the center, which means inserting a triangle in place of the center point (see Figure \ref{blowup2}, first line). 
Though this is not so honest as 3-dimensional pictures, this is still a fairly honest way to
describe the cutting process. 
\par\noindent
There is a way to make these pictures simpler without loosing any information: 
We watch only one chamber, say 3, and show the cutting/blowing-up of the center only in the chamber 3 (see Figure \ref{blowup2}, second line, right). You see that the terminal figure of marked triangle tells everything.
\begin{figure}[h]
\begin{center} 
\include{blowup2} 
\end{center}
\caption{Blowing up a point to be a triangle (blowup2)}
\label{blowup2}
\end{figure} 
\par\smallskip
Now we are ready. Label the tetrahedra as $1,2,\dots,5$ as in the previous subsection. We choose the edge $E$ joining two vertices 1235 and 1345. Along this edge $E$, 
there are three tetrahedra $1, 3$ and 5. Other than these three, tetrahedron 2 touches the vertex 1235, and tetrahedron 4 touches the vertex 1345. 
We cut the edge by a new hyperplane 6. The tetrahedron 2 (resp. 4) has vertex 1235 (resp. 1345) as the intersection with $E$; 
this versex is cut, and the tetrahedron 2 becomes a prism (see Figure \ref{tetratoprism1}).
\begin{figure}[h]
\begin{center}
\include{tetratoprism1} 
\end{center}
\caption{Truncating a vertex of a tetrahedron (tetratoprism1)}
\label{tetratoprism1}
\end{figure} 
The tetrahedron 5 (as well as 1 and 3) has the edge $E$; this edge is cut, and the tetrahedron 5 becomes a prism (see Figures \ref{tetratoprism2} and \ref{prismcut}).
\begin{figure}[h]
\begin{center} 
\include{tetratoprism2} 
\end{center}
\caption{Truncating an edge of a tetrahedron (tetratoprism2)}
\label{tetratoprism2}
\end{figure} 
The new face 6 is a prism bounded by two triangles labeled 2 and 4, and three rectangles labeled 1, 3 and 5. Recall the above 1-dimension lower process, then the blow-up of the edge $E$ to become the prism 6 can be summerized by Figure \ref{blowup56}.
\begin{figure}[h]
\begin{center} 
\include{blowup56} 
\end{center}
\caption{Blowing up a segment to be a prism (blowup56)}
\label{blowup56}
\end{figure} 
In the case of seven or more hyperplanes, we describe the cutting/blow-up process only by this kind of figure.
\subsection{Seven hyperplanes in the 4-space}
Seven hyperplanes in the projective 4-space cut out a {\it unique} chamber CC stable under the action of $\mathbb{Z}_7$. It is bounded by the seven dumplings. We study how the seven dumplings tessellate a 3-sphere.
\subsubsection{Labeling and a study of 2-skeleton}
The central chamber CC has seven faces $D_1,\dots,D_7$, which form the $\mathbb{Z}_7$-orbit of the dumpling
\small$$D_1={\rm CC}\cap H_1:\quad\scriptsize\begin{array}{ccccccccc}
    & &    & &1457&   &    &   &\\          
    & &/\quad   &4 &|   & 5  &\quad\ba &   &\\          
1347& &   & &1245&   &     &   &1567\\
    &\quad\ba&    &/&    &\ba& &/\quad&\\ 
\quad \ba   & 3  &1234& &    &   &1256& 6 &/\quad\\
    &\ba\quad&| & &  2 &   & |  & \quad/ &\\ 
    & &1237&-&--  & -&1267&   &      
\end{array}.$$\normalsize
Note that the pentagonal face of $D_1$ behind is $D_1\cap D_7$ and the pentagon in front is $D_1\cap D_2$. 
\footnote{Using the notation introduced in \cite{CY}, CC is represented by $+++++++$, and its boundary consistes of 
$$D_1=0++++++,\quad D_2=-0++++,\ \dots,\ D7=------0.$$
The boundary of $D_1$ consists of
$$\begin{array}{lll}
00+++++:{\rm5-gon}&
0-0++++:{\rm triangle}&
0--0+++:{\rm rectangle}\\
0---0++:{\rm rectangle}&
0----0+:{\rm triangle}&
0-----0:{\rm5-gon}
\end{array}
$$} 
The 2-skeleton of CC consists of seven pentagons, seven triangles and seven rectangles.  Remove the intersections  of the two consecutive ones:
$$D_1\cap D_2,\quad D_2\cap D_3,\quad \dots,\quad D_7\cap D_1:$$
\comment{
$$\scriptsize\begin{array}{cccccccccc}
         &                &2367-2356          &                &4512-4571&                &6734-6723&                &\\  
         &\qquad\qquad/   &                   &\ba\qquad\qquad/&         &\ba\qquad\qquad/&         &\ba\qquad\qquad &  \\
         &7123            & 23                &2345            & 45      &4567            &  67     &6712            &\\
         &/\qquad\qquad\ba&                   &/\qquad\qquad\ba&         &/\qquad\qquad\ba&         &/\qquad\qquad\ba&\\         
6712     &  12            &1234               &  34            &3456     &  56            &5671     &71              &7123   \\
\qquad\ba&                &/\qquad\qquad\ba&  &/\qquad\qquad\ba&         &/\qquad\qquad\ba&         &/\qquad\qquad   &\\
         &1256-1245       &                   &3471-3467       &         &5623-5612       &         &1457-1347       &
\end{array}\normalsize$$}
$$\scriptsize\begin{array}{cccccccccc}
         &                &2367-2356          &                &\cdots          &                &\\  
         &\qquad\qquad/   &                   &\ba\qquad\qquad/&                &\ba\qquad\qquad &  \\
         &7123            & 23                &2345            & \cdots         &6712            &\\
         &/\qquad\qquad\ba&                   &/\qquad\qquad\ba&                &/\qquad\qquad\ba&\\         
6712     &  12            &1234               &  34            &\cdots          &71              &7123   \\
\qquad\ba&                &/\qquad\qquad\ba   &                &/\qquad         &\ba\qquad\qquad\qquad/&  \\
         &1256-1245       &                   &3471-3467       &                &1457-1347&
\end{array}\normalsize$$
\normalsize
(Note that the union of these seven pentagons form a M\"obius strip, but for a while just forget it.)
Then  the remaining  seven triangles and seven rectangles form a M\"obius strip
$$M:\quad\scriptsize\begin{array}{ccccccccccc  }
6712&-&7123&-&1234&-&2345&-&3456&-&4567\\
 \quad\ba&27 & /\quad\ba &31&/\quad\ba  &42&/\quad\ba&53&/\quad\ba&64&/\quad \\
    &2367&&3471& &1245& &2356& &3467&\\
 &\quad\ba&37 & /\quad\ba &41&/\quad\ba  &52&/\quad\ba   &63&/\quad   & \\
 &   &6734&&7145& &5612& &2367&&\\
 &&\quad\ba&47 & /\ 57\ \ba &51&/\ 61\ \ba  &62&/\quad   &&    \\
&&&4567&-&5671&-&6712&&&
\end{array}\normalsize$$
\begin{center} M\"obius  7
\end{center}
Its boundary 
$$C=\partial M:\quad 1234-2345-3456-4567-5671-6712-7123-1234$$
is an unknotted circle because $C$ bounds also the disc
$$\scriptsize\begin{array}{cccccccc}
7123&-&-&1234&-&-&2345\\[2mm]
|&\ba&13&/\qquad\ba&24&/&|\\[2mm]
|&&3471&14&1245&&|\\[2mm]
|&&&\ba\qquad/&&&|\\[2mm]
6712&&71&7145&45&&3456\\[2mm]
&\ba&&/\quad57\quad\ba&&/&\\
&&5671&-&4567&&&
\end{array}\normalsize.$$

\subsubsection{From a solid torus}
We start with a vertical pentagonal prism, horizontally sliced into seven thinner pentagonal prisms. 
Its top and the bottom pentagons are identified after a $4\pi/5$ rotation; 
this makes the vertical prism a solid torus with a tessellation of the lateral boundary into 35 rectangles. Now collapse vertically seven diagonal rectangles. 
This changes 21 rectangles (out of 35) into seven segments (forming a circle) and fourteen triangles. 
Accordingly, the seven thinner pentagonal prisms  into dumplings. Next use each collapsed edge as a valley of each pair of triangles 
face-to-face in order to close the pair like a book, and further identify seven pairs of rectangles. 
This identification changes the solid torus into a sphere tessellated by seven dumplings.
\par\medskip\noindent{\bf A vertical pentagonal prism horizontally sliced.}
We show the above process by using the labels $1,2,3,4,5,6,7$ of the seven hyperplanes. The boundary of the vertical pentagonal prism; 
each line, say $71$ is a pentagon on the plane $7\cap1$. Two pentagons $71$ and $12$ sandwich a thin pentagonal prism. $6712$ is the point $6\cap7\cap1\cap2$.

$$\scriptsize\begin{array}{rccccccccccc}
71:&7126&-&7123&-&7134&-&7145&-&7156&-&7126\\
  &\|   &\times&\|  &13&|   &14&|   &15&|   &16&\|\\
12:&1267&-&1237&-&1234&-&1245&-&1256&-&1267\\
  &|   &27  &\| &\times&\|   &24&|   &25&|   &26&|\\
23:&2367&-&2371&-&3412&-&2345&-&2356&-&2367\\
  &|   &37  &|  &31&\| &\times&\|   &35&|   &36&|\\
34:&3467&-&3471&-&3412&-&3452&-&3456&-&3467\\
  &|   &47  &|  &41&| &42&\|   &\times&\|   &46&|\\
45:&4567&-&4571&-&4512&-&4523&-&4563&-&4567\\
  &\|   &57  &|  &51&| &52&|   &53&\|   &\times&\|\\
56:&5674&-&5671&-&5612&-&5623&-&5634&-&5674\\
  &\|   &\times&\|  &61&| &62  &|   &63&|   &64&\|\\
67:&6745&-&6715&-&6712&-&6723&-&6734&-&6745\\
  &|   &75&\|  &\times&\| &72  &|   &73&|   &74&|\\
71:&7145&-&7156&-&7126&-&7123&-&7134&-&7145\\
  &|   &15  &|  &16&\| &\times  &\|   &13&|   &14&|\\
\end{array}\normalsize$$
\begin{center} Table of piles 7
\end{center}
\par\medskip\noindent{\bf Collapsing rectangles.}
If we identify two points with the same {\it consecutive} indices {\it as sets}, for exmaple, $7123=1237=2371$ but $2356\not=5623$,
you find seven rectangles with pairwise coincide vertices (marked by $\times$), situated diagonally. We collapse such rectangles vertically to a segment, each called a {\bf special edge}. Note that a special edge has two ends which are labeled by consecutive numbers: $\{1234,2345\},\{2345,3451\},\dots$ Accordingly the remaining rectangles next to the special edges become triangles, and the seven thin rectangular prisms become dumplings. But as a whole it still remains to be a cylinder.
\par\medskip\noindent{\bf Identifying the top and the bottom to get a solid torus $ST$.} 
We identify the top pentagon and the bottom pentagon of the collapsed vertical prism according to the labels; note that we must twist the triangle by $4\pi/5$. By the identification we get a solid torus, say $ST$; its meridean ${\bf mer}$ is represented by the boundary of a(ny) pentagon, in Table of piles 7, say 71 (or 12, 23,\dots). Note that it is still an abstract object. The seven special edges now form a circle:
$$C: 1234-2345-3456-4561-5612-6712-7123-1234.$$
We take a vertical line, in Table of piles 7, with a $4\pi/5$ twist, say,
$${\bf par}: 7123 - 3471 - 4571 - 5671 - 7123$$
as a parallel.

\par\medskip\noindent{\bf The torus $T$ and the solid torus $UT$.}
We consider in our space ($\sim \mathbb{S}^3$) a usual (unknotted) solid torus, say $UT$, and think our solid torus $ST$ fills outside of $UT$. Let us see the boundary torus
$$T:=\partial\ ST\ (=\partial\ UT)$$
has {\bf par} as the meridean, and {\bf mer} as the parallel.
From the   description of  the vertical prism, we see fourteen triangles and fourteen rectangles tessellating the torus $T$ arranged in a hexagonal way (identify the opposite sides of the hexagon):
$$\scriptsize\begin{array}{ccccccccccc}
6712&-&6723&-&6734&-&4567& & &&\\
|   &27\ba72&|  &73&|   &74&|   &75\ba&   &&\\
2367&-&7123&-&7134&-&7145&-&5671&&\\
|   &37 &| &31\ba13&|   & 14&|   &15 &|   &16\ba&\\
3467&-&3471&-&1234&-&1245&-&1256&-&6712\\
|   &47 &|  &41 &| &42\ba24&|   &25 &|   & 26&|\\
4567&-&4571&-&4512&-&2345&-&2356&-&2367\\
    &\ba57 &|  &51 &| &52 &|   &53\ba35&|   &36 &|\\
    & &5671&-&5612&-&5623&-&3456&-&3467\\
    &   &   &\ba61&| &62 &|   &63 &|   &64\ba46&|\\
    &   &   &   &6712&-&6723&-&6734&-&4567\\
\end{array}\normalsize$$
\begin{center} Torus (hexagon) 7
\end{center}

\par\medskip\noindent{\bf Collapsing the solid torus $UT$ by folding $T$.} We identify two vertices with the same indices as sets. Along each collapsed edge (special edge), there are two triangles with the same vertices; these triangles are folded along the special edge into one triangle. Further there are seven pairs of rectangles with the same vertices; they are also identified. Consequently the torus $T$ is folded along the curve $C$ to be a M\"obius strip, which turns out to be the same as the one in the previous subsection (shown as M\"obius 7).
\comment{$$M:\quad\scriptsize\begin{array}{ccccccccccc  }
6712&-&7123&-&1234&-&2345&-&3456&-&4567\\
 \quad\ba&27 & /\quad\ba &31&/\quad\ba  &42&/\quad\ba&53&/\quad\ba&64&/\quad \\
    &2367&&3471& &1245& &2356& &3467&\\
 &\quad\ba&37 & /\quad\ba &41&/\quad\ba  &52&/\quad\ba   &63&/\quad   & \\
 &   &6734&&7145& &5612& &2367&&\\
 &&\quad\ba&47 & /\ 57\ \ba &51&/\ 61\ \ba  &62&/\quad   &&    \\
&&&4567&-&5671&-&6712&&&
\end{array}\normalsize$$
\begin{center} M\"obius  7
\end{center}}
In this way, we get a space $X_7$ from the solid torus $ST$ made by six prisms by folding the torus $T$ (collapsing $UT$).

Since the curve $C$ is unknotted, this folding is done exactly the same as in the previous section. Then $X_7$ is homeomorphic to $\mathbb{S}^3$.

\par\bigskip\noindent{\bf A model.} In \cite{CY}, the tessellation is shown, as we quote in Figure \ref{7facebody}. In the figure, 0 stands for 7, and the edge connecting $0346$ and $0134$ is not drawn to avoid complication. Six dumplings are shown; they are packed in a dumpling. The other one is infinitely large. From this diagram, one can hardly see the $\mathbb{Z}_7$ action. To understand this group action is the motivation of our prism-collapse construction.
\begin{figure}[t]
\begin{center}
\include{7facebody70}
\end{center}
\caption{Six dumplings pack a dumpling}
\label{7facebody}
\end{figure}

\subsubsection{$7-1=6$}
If we remove the seventh hyperplane, then the remaining six hyperplanes bounds a hyperhexahedron. This process can be described as follows: the dumpling on the seventh plane reduces to the union of two adjacent segments. More precisely, the two triangles of the boundary of the dumpling reduce to two points, the segment joining the remaining two vertices to a point, and the dumpling to the union of two segments joining the last one to the former ones. Combinatorial explanation: The dumpling on the seventh plane with vertices
\small$$\scriptsize\begin{array}{ccccc}
&&3467&&\\
&/&|&\ba&\\
2367&&|&&4567\\
|&\ba\qquad&|&\qquad/&|\\
|&1267&---&1567&|\\
|&/\qquad&|&\qquad\ba&|\\
1237&&|&&1457\\
&\ba&|&/&\\
&&1347\end{array}\normalsize\quad {\rm on\ 7}$$\normalsize
reduces to the union of two segments 
$$\bf1236\ --\ 1346\ --\ 1456.$$
Note that 
$$\begin{array}{ll}
\{1,2,3,6\}&=\{2,3,6,7\}\cup\{1,2,6,7\}\cup\{1,2,3,7\}-\{7\},\\[2mm]
\{1,3,4,6\}&=\{1,3,4,7\}\cup\{3,4,6,7\}-\{7\},\\[2mm]
\{1,4,5,6\}&=\{4,5,6,7\}\cup\{1,5,6,7\}\cup\{1,4,5,7\}-\{7\}.
\end{array}$$
Accordingly the other dumplings reduce to prisms (cf. Figure \ref{prismcut}), for example,
\small
$$
\scriptsize\begin{array}{ccc}
&1245&\\
&/\qquad|\qquad\ba&\\
7145&|&2345\\
|&\ba\qquad|\qquad/&|\\
|&6745--6345&|\\
|&/\qquad|\qquad\ba&|\\
6715&|&6235\\
&\ba\qquad|\qquad/&\\
&6125&\end{array}\normalsize
\longrightarrow
\scriptsize\begin{array}{ccc}
&1245&\\
&/\qquad|\qquad\ba&\\
\quad\qquad/&|&2345\\
\qquad/&\qquad|\qquad/&|\\
\bf1465&----6345&|\\
\qquad\ba&\qquad|\qquad\ba&|\\
\quad\qquad\ba&|&6235\\
&\ba\qquad|\qquad/&\\
&6125&\end{array}\normalsize
\quad
\normalsize{\rm on\ 5} \ ({\rm similar\ on}\ 2),
$$
$$
\scriptsize\begin{array}{ccc}
&6723&\\
&/\qquad|\qquad\ba&\\
5623&|&7123\\
|&\ba\qquad|\qquad/&|\\
|&4523--4123&|\\
|&/\qquad|\qquad\ba&|\\
4563&|&4713\\
&\ba\qquad|\qquad/&\\
&4673&\end{array}\normalsize
\longrightarrow
\scriptsize\begin{array}{ccc}
5623&----&\bf1236\\
|&\ba\qquad\qquad/&|\\
|&4523--4123&|\\
|&/\qquad\qquad\ba&|\\
4563&----&\bf1346\\
\end{array}\normalsize
\quad
{\rm on\ 3}\ ({\rm similar\ on}\ 4),
$$
$$
\scriptsize\begin{array}{ccc}
&2356&\\
&/\qquad|\qquad\ba&\\
1256&|&3456\\
|&\ba\qquad|\qquad/&|\\
|&7156--7456&|\\
|&/\qquad|\qquad\ba&|\\
7126&|&7346\\
&\ba\qquad|\qquad/&\\
&7236&\end{array}\normalsize
\longrightarrow
\scriptsize\begin{array}{ccccc}
&&2356&&\\
&/&/&\ba&\\
1256&&&&3456\\
|&\ba&&/&\\
|&/\quad&\bf1456&\quad/&\\
\bf1236&&|&&\\
&\ba&|&/&\\
&&\bf1346&&\end{array}\normalsize
\quad{\rm on\ 6}\ ({\rm similar\ on}\ 1).$$
\normalsize
\subsubsection{$6+1=7$\label{67}}
Let us describe the change of 3-sphere tessellated by six prisms into that by seven dumplings. As we explained in \S\ref{56}, 
we express this change -- blowing-up a connected union of two edges (of a prism) -- by cutting a prism 6 by a hyperplane 7 and dividing the prism into two dumplings, 
which share the pentagonal face that should be called 6. (See Figure \ref{blowup67}).
\begin{figure}[h]
\begin{center} 
\include{blowup67} 
\end{center}
\caption{Blowing up the union of two segments to be a dumpling (blowup67)}
\label{blowup67}
\end{figure} 

\subsection{$m(\ge8)$ hyperplanes in the 4-space}
A Veronese arrangement of $m(\ge8)$ hyperplanes in the projective 4-space cut out a {\it unique} chamber CC, stable under the action of $\mathbb{Z}_m$. 
It is bounded by $m$ $(m-2)$-dumplings. (An $(m-2)$-dumpling is bounded by two $(m-2)$-gons, two triangles and $m-5$ rectangles.) 
The reader is expected to imagine how these $(m-2)$-dumplings tessellate a 3-sphere. Here we show a 6-dumpling ($m=8$):
\small
$$D_1={\rm CC}\cap H_1:\quad\scriptsize\begin{array}{ccccccc}
 1458   &- &  -   &-&    -   & -  &1568\\          
  |  & \ba&    & 5  &    & /  &|\\          
|& & 1245   &--&      1256  &   &|\\
    |&4&  |  &    & |&6&|\\ 
1348   & -  &1234 &  2  &   1267& - &1678\\
    &\ba\quad3&|  &   &   |  &7\quad/ &\\ 
    & &1238&--  & 1278&   &      
\end{array}\normalsize.$$
\normalsize

The central chamber CC has $m$ faces $D_1,\dots,D_m$, on which $\mathbb{Z}_m$ acts as $D_k\to D_{k+1}$ mod $m$. The 2-skeleton of CC consists of $m$ $(m-2)$-gons, $m$ triangles and $m(m-5)/2$ rectangles. Remove the intersections of the two consecutive ones:
$$D_1\cap D_2,\quad D_2\cap D_3,\quad \dots,\quad D_m\cap D_1.$$
Then the remaining $m$ triangles and $m(m-5)/2$ rectangles form a M\"obius strip $M$, which the reader can make following the previous sections. 
\par\noindent
Its boundary
$$C=\partial M:\quad 1234-2345-\cdots-m123-1234$$
is an unknotted circle because $C$ bounds also a disc. To find a disc is not so obvious so we give a hint when $m=8$ (the reader is expected to guess the codes for $A,B,C,D,E$, and then to guess what happens when $m\ge9$)
\small$$\scriptsize\begin{array}{cccccccc}
8123&-&1234&-&2345&-&3456\\[2mm]
|&\ba&/\quad\quad\ba&&/\quad\quad\ba&/&|\\[2mm]
|&\qquad A&&B&&C\qquad&|\\[2mm]
|&&\ba\quad\quad/&&\ba\quad\quad/&&|\\[2mm]
|&&D&15&E&&|\\[2mm]
|&&&\ba\qquad/&&&|\\[2mm]
7812&&81&8156&56&&4567\\[2mm]
&\ba&&/\quad68\quad\ba&&/&\\
&&6781&-&5678&&&
\end{array}\normalsize.$$
(Answer: $A=3481,B=1245,C=2356,D=8145,E=1256$.)
\normalsize
\subsubsection{From a solid torus}
You start with a vertical $(m-2)$-gonal prism, horizontally sliced into $m$ thinner $(m-2)$-gonal prisms. Its top and the bottom triangles are identified after a $4\pi/(m-2)$ rotation. This makes the vertical prism a solid torus; the boundary is tessellated by $m(m-2)$ rectangles. 
\par\smallskip\noindent
Now collapse vertically $m$ diagonal rectangles. This changes $3m$ rectangles into $m$ segments (special edges) and $2m$ triangles. Accordingly, the $m$ thinner pentagonal prisms deform into $(m-2)$-dumplings, forming a solid torus $ST$. The special edges make a circle $C$.
\par\smallskip\noindent
We consider in our space ($\sim \mathbb{S}^3$) a usual solid torus, say $UT$, and think our solid torus $ST$ fills outside of $UT$.
\par\smallskip\noindent
In the tessellation on the boundary torus
$$T:=\partial\ ST= \partial\ UT,$$
each special edge is shared by the two triangles with the same code. The special edge is used as a valley to close the pair of triangles like a book, and further identify  pairs of the rectangles. This identification changes the torus $T$ into the M\"obius strip $M$ with boundary $C$, and accordingly collapse the solid torus $UT$, and changes the solid torus $ST$ into a 3-sphere tessellated by $m$ $(m-2)$-dumplings, recovering the tessellation on $\partial\ {\rm CC}$.

\subsubsection{From $m$ hyperplanes to $m-1$ hyperplanes}
Consider the central chamber bounded by $m$ hyperplanes. If we remove the $m$-th hyperplane, 
then $m-1$ hyperplanes bound the central chamber bounded by $m-1$ hyperplanes. This process can be described as follows: 
the $(m-2)$-dumpling on the $m$-th plane reduces to $m-4$ connected segments.  Combinatorial explanation: The $(m-2)$-dumpling on the $m$-th plane with vertices
$$\scriptsize\begin{array}{ccccccc}
&&\bullet&-\ \bullet\ \cdots\ -\cdots\ \bullet\ -&\bullet&&\\
\bullet&&|&&|&&\bullet\\
|&\bullet&-&--\cdots\ -\cdots\ --&-&\bullet&|\\
\bullet&&|&&|&&\bullet\\
&&\bullet&-\ \bullet\ \cdots\ -\cdots\ \bullet\ -&\bullet&&
\end{array}\normalsize$$
reduces to the union of $m-4$ segments 
$$\bullet \ -\ \ \bullet\ \ -\ \ \bullet\ \cdots\ -\cdots\ \ \bullet\ -\ \ \bullet\ -\ \ \bullet.$$
Accordingly, the other $(m-2)$-dumplings reduce to $(m-3)$-dumplings.

\subsubsection{From $m$ hyperplanes to $m+1$ hyperplanes}
Let us describe the change of 3-sphere tessellated by $m$ $(m-2)$-dumplings into that by $m+1$ $(m-1)$-dumplings. 
As we explained in \S\ref{56}, we express this change -- blowing-up a connected union of $m-4$ edges (of a $(m-2)$-dumpling) -- by cutting the $(m-2)$-dumpling 
(labeled by $m$) by the 
hyperplane (labeled by $m+1$) and dividing the $(m-2)$-dumpling into two $(m-1)$-dumplings, which share the $(m-1)$-gonal face which should be labeled by $m$. 
(See Figure \ref{blowup89}).
\begin{figure}[h] 
\begin{center} 
\include{blowup89} 
\end{center}
\caption{Blowing up the union of $m-4$ segments to be a $(m-1)$-dumpling (blowup89)}
\label{blowup89}
\end{figure}  
\clearpage
\section{Appendix: Higher dimensional cases}
We describe higher dimensional dumplings $D_n$ for odd $n$, and central chambers CC$_n$ for even $n$ in arrangements of $n+3$ hyperplanes 
in $\mathbb{P}^n$. Description will be {\bf just combinatorial}. We start from dimension 0: 0-dimensional central chamber CC$_0$ is a point.
\subsection{Review}
{\bf - D$_1$:} 1-dimensional dumpling is a segment, whose boundary consists of two CC$_0$'s.

\noindent
{\bf - CC$_2$:} The 2-dimensional central chamber CC$_2$ (usually called a pentagon) is bounded by five $D_1$'s. Each $D_1$ share CC$_0$-faces with two $D_1$'s.

\noindent
{\bf - D$_3$:}
Consider the direct product of CC$_2$ and an interval $I$. Choose a boundary component of CC$_2$, and call it ${\bf D_1}$. 
Push down ${\rm\bf D_1}\times I $ to ${\rm\bf D_1}$, to get a dumpling $D_3$. The crashed ${\bf D_1}$ is called the {\bf special} edge of $D_3$.  Accordingly, for a boundary component $D_1$ adjacent to ${\bf D_1}$ (there are two of them), the rectangle ${ D_1}\times I $ is pushed to become a triangle.  For a boundary component $D_1$ non-adjacent to ${\bf D_1}$ (there are two of them), the rectangle ${ D_1}\times I $ leaves as it is.
Thus the boundary of a dumpling consists of two CC$_2$'s and two triangles and two quadrilaterals (The following notation is explained in \cite{CY}):
$$\begin{array}{lll}
\partial(++++++)&=0+++++ &\cdots\ \rm CC_2\\
 &\cup\ -0++++ &\cdots\ \rm triangle \\ 
&\cup\ --0+++ &\cdots\ \rm quadrilateral\\
& \cup\ ---0++ &\cdots\ \rm quadrilateral\\
&\cup\ ----0+ &\cdots\ \rm triangle\\
& \cup\ -----0 &\cdots\ \rm CC_2 
\end{array}\normalsize$$

\noindent
{\bf - CC$_4$:} \label{CC4} 
The 4-dimensional central chamber CC$_4$ is bounded by seven $D_3$'s. Each $D_3$ shares CC$_2$-faces with two $D_3$'s. The 2-skeleton of the $D_3$-tessellation of $\mathbb{S}^3$ minus the seven CC$_2$'s form a M\"obius strip $M_2$ with boundary $C_1$ consisting of seven special edges $D_1$'s of seven $D_3$'s. 
\par\medskip\noindent
On the other hand, we start with a vertical prism with base CC$_2$, horizontally sliced into seven thinner prisms with base CC$_2$. Collapse the seven diagonal $D_1\times I$'s on the boundary into seven $\bf D_1$'s; 
this makes each thinner prism a $D_3$ with the said $\bf D_1$ as special one. Identify the top and the bottom CC$_2$ after $4\pi/5$ rotation so that the seven special edges $\bf D_1$ form a circle $C_1$. In this way, we get a solid torus $ST_3$ made of seven $D_3$'s. 
\par\medskip\noindent
The (tessellated) boundary torus $T_2$ of $ST_3$ is just the double of the (tessellated) M\"obius strip $M_2$ branching along $C_1$.
\par\medskip\noindent
Fold the boundary torus $T_2$ along the curve $C_1$ to $M_2$, Then the solid torus $ST_3$ changes into the tessellated boundary of CC$_4$.

\subsection{$D_{2k-1},D_{2k}$ $(k\geqslant3)$}

\noindent
{\bf - D$_{2k-1}$:}
Consider the direct product of CC$_{2k-2}$ and an interval $I$. Choose a boundary component $D_{2k-3}$ of CC$_{2k-2}$, and call it ${\bf D_{2k-3}}$. 
Push down ${\bf D_{2k-3}}\times I $ to ${\bf D_{2k-3}}$, to get a dumpling $D_{2k-1}$. 
The resulting {$\bf D_{2k-3}$} is called the {\bf special} edge of $\bf D_{2k-1}$.  

\noindent
{\bf - CC$_{2k}$:}
The $2k$-dimensional central chamber CC$_{2k}$ is bounded by $2k+3$ dumplings $D_{2k-1}$. 
Each $D_{2k-1}$ shares CC$_{2k-2}$-faces with two $D_{2k-1}$'s. 
The ($2k-2$)-skeleton of the $D_{2k-1}$-tessellation of $\mathbb{S}^{2k-1}$ minus the $2k+3$ chambers CC$_{2k-2}$ form a 
CW-complex $M_{2k-2}$ with boundary $C_{2k-3}$ $(\sim \Delta_{2k-4}\times \mathbb{S}^1$) consisting of $2k+3$ special edges $\bf D_{2k-3}$ of $2k+3$ dumplings $D_{2k-1}$. 
\par\medskip\noindent
On the other hand, we start with a vertical prism with base CC$_{2k-2}$, horizontally sliced into $2k+3$ thinner prisms with base CC$_{2k-2}$. 
Collapse the $2k+3$ diagonal $D_{2k-3}\times I$ on the boundary into $2k+3$ dumplings $\bf D_{2k-3}$; 
this makes each thinner prism a dumpling $D_{2k-1}$ with the said $\bf D_{2k-3}$ as { special} edge. 
Identify the top and the bottom CC$_{2k-2}$ after $4\pi/(2k+1)$ rotation so that the $2k+3$ special edges $\bf D_{2k-3}$ form $C_{2k-3}$. 
In this way, we get a solid torus $ST_{2k-1}\ (\sim \Delta_{2k-2}\times \mathbb{S}^1)$ made of $2k+3$ dumplings $D_{2k-1}$. 
\par\medskip\noindent
The (tessellated) boundary torus $T_{2k-2}\ (\sim \mathbb{S}^{2k-3}\times \mathbb{S}^1)$ of $ST_{2k-1}$ can be folded  
along the solid ($2k-3$)-torus $C_{2k-3}$ (as a hinge)
to $M_{2k-2}$. Then the solid torus $ST_{2k-1}$ changes into the tessellated boundary of CC$_{2k}$.

\par\bigskip\noindent
{\bf Acknowledgement:} The last author is grateful to K. Cho for his help. This collaboration started on 11 March 2011
, the very day of the Fukushima disaster in Japan.

\medskip
\begin{flushleft}
Francois Ap\'ery\\
Universit\'e de Haute-Alsace\\
2, rue des Fr\`eres Lumi\`ere\\
68093 Mulhouse Cedex\\
France\\
francois.apery@uha.fr
\par\medskip\noindent
Bernard Morin\\
villa Beausoleil\\
32, avenue de la R\'esistance\\
92370 Chaville\\
France\\
\par\medskip\noindent
Masaaki Yoshida\\
Kyushu University\\
Nishi-ku, Fukuoka 819-0395\\
Japan\\ 
myoshida@math.kyushu-u.ac.jp
\end{flushleft}

\end{document}